\documentclass[10pt]{article}

\usepackage[latin1]{inputenc}
\usepackage{tikz}
\usetikzlibrary{arrows.meta,calc}
\usepackage[pagebackref=true, colorlinks=true, linkcolor=blue, citecolor=blue, linkcolor=blue, anchorcolor=red, urlcolor=blue]{hyperref}
\usepackage{url, amsmath,enumerate,fancyhdr,amssymb, amsthm, 
pstricks, epsf, lscape, ifpdf, graphicx, float}
\usepackage{algorithm}
\usepackage{enumitem}
\usepackage{color}
\usepackage{algpseudocode} 
\algnewcommand\algorithmicinput{\textbf{Input:}}
\algnewcommand\INPUT{\item[\algorithmicinput]}
\algnewcommand\algorithmicinitialization{\textbf{Initialization:}}
\algnewcommand\INITIALIZATION{\item[\algorithmicinitialization]}

\newcommand{\ti}{\times}

	\newcommand\bovermat[2]{%
		\makebox[0pt][l]{$\smash{\overbrace{\phantom{%
						\begin{matrix}#2\end{matrix}}}^{\text{#1}}}$}#2}
	
	\usepackage{mathtools}
	\makeatletter
	\renewcommand*\env@matrix[1][*\c@MaxMatrixCols c]{%
		\hskip -\arraycolsep
		\let\@ifnextchar\new@ifnextchar
		\array{#1}}
	\makeatother
	
	\usepackage{tikz,array,calc, xcolor}
	\usetikzlibrary{decorations.pathreplacing}

\usepackage[margin=3cm]{geometry}

\mathchardef\mhyphen="2D

\newtheorem{Definition}{Definition}

\newtheorem{Example}{Example}
\newtheorem{Proposition}{Proposition}
\newtheorem{Lemma}{Lemma}
\newtheorem{Theorem}{Theorem}
\newtheorem{Corollary}{Corollary}
\newtheorem{Remark}{Remark}
\newtheorem{Assumption}{Assumption}

\newcommand{\ba}{\begin{array}}  
	\newcommand{\ena}{\end{array}}

\newcommand{\A}{\mathcal A}

\newcommand{\bpx}{\begin{pmatrix}}
\newcommand{\epx}{\end{pmatrix}}
\newcommand{\bbx}{\begin{bmatrix}}
\newcommand{\ebx}{\end{bmatrix}}

\newcommand{\bdef}{\begin{Definition}} 
\newcommand{\commentout}[1]{}
\newcommand{\co}[1]{}

\newcommand{\coab}[1]{}

\newcommand{\norm}[1]{\parallel \! #1 \! \parallel}
\newcommand{\sym}[1]{{\cal S}^{#1}}
\newcommand{\psd}[1]{{\cal S}_+^{#1}}
\newcommand{\pd}[1]{{\cal S}_{++}^{#1}}
\newcommand{\mymatrix}[1]{{\cal R}_+^{#1}}
\newcommand{\rad}[1]{\mathbb{R}^{#1}}

\newcommand{\symn}{\sym{n}}
\newcommand{\psdn}{\psd{n}}

\newcommand{\N}{ {\cal N} }

\newcommand{\J}{\mathcal J}

\newcommand{\la}{\langle}
\newcommand{\ra}{\rangle}

\newcommand{\bY}{\bar{Y}}

\newcommand{\dist}{\operatorname{dist}}

\newcommand{\beq}{\begin{equation}}
\newcommand{\eeq}{\end{equation}}
\newcommand{\beqa}{\begin{eqnarray}}
\newcommand{\eeqa}{\end{eqnarray}}
\newcommand{\bac}{\begin{array}{ccccccccccc}}
\newcommand{\eac}{\end{array}}
\newcommand{\bprop}{\begin{Proposition}}
\newcommand{\eprop}{\end{Proposition}}

\newcommand{\beqast}{\begin{eqnarray*}}
\newcommand{\eeqast}{\end{eqnarray*}}
\newcommand{\benum}{\begin{enumerate}}
\newcommand{\eenum}{\end{enumerate}}
\newcommand{\bit}{\begin{itemize}}
\newcommand{\eit}{\end{itemize}}
\newcommand{\bth}{\begin{Theorem}}
\newcommand{\enth}{\end{Theorem}}
\newcommand{\ble}{\begin{Lemma}}
\newcommand{\ele}{\end{Lemma}}
\newcommand{\bex}{\begin{Example}}
\newcommand{\eex}{\end{Example}}
\newcommand{\bcor}{\begin{Corollary}}
\newcommand{\ecor}{\end{Corollary}}
\newcommand{\brem}{\begin{Remark}}
\newcommand{\erem}{\end{Remark}}
\newcommand{\bass}{\begin{Assumption}}
\newcommand{\eass}{\end{Assumption}}

\newcommand{\val}{\operatorname{val}}

\newcommand{\bsmx}{\begin{small} \begin{pmatrix}}
\newcommand{\esmx}{\end{pmatrix} \end{small}}

\title{\bf A combinatorial approach to Ramana's exact dual for semidefinite programming}

\author{G\'{a}bor  Pataki \thanks{Department of Statistics and Operations Research, University of North Carolina   at Chapel Hill} \hspace{1cm} 
}

\begin{document}

\maketitle

\begin{abstract}
Thirty years ago, in a seminal paper Ramana derived an exact dual for Semidefinite Programming (SDP). Ramana's dual has the following remarkable features:
i) it is an explicit, polynomial size  semidefinite program ii) it does not assume that the primal is strictly feasible, nor does it make any other regularity assumptions iii) yet, it has strong duality with the primal. The complexity implications of  Ramana's dual are  fundamental, 
and to date still the best known. The most important of these is that SDP feasibility in the Turing model is not NP-complete, unless NP = co-NP.

We give a  treatment of Ramana's dual which is both 
simpler and more complete,  than was previously available.
First we connect it to a seemingly very different way of inducing strong duality: 
reformulating the  SDP into a rank revealing form  using  elementary row operations
and rotations.
Second,  while previous works characterized its objective value, 
we completely characterize its feasible set: in particular, we show 
it is a higher dimensional representation of 
an exact dual, which, however is not an explicit SDP. 
We also prove that -- somewhat surprisingly -- strict feasibility of Ramana's dual implies that the only feasible solution of the primal is the zero matrix.

As a corollary, we obtain a short and transparent derivation of Ramana's dual, which we believe  is accessible to both the optimization 
and the theoretical computer science communities.  
Our approach  is combinatorial in the following sense: i) we use a minimum amount of continuous optimization theory ii) 
we show that feasible solutions in Ramana's dual are identified with regular facial reduction sequences, i.e., essentially discrete structures.

\vspace{.2cm} 

\vspace{.2cm} 
\end{abstract}

\vspace{.2cm} 
{\bf Dedicated to Motakuri Venkata Ramana}
\tableofcontents

\section{Introduction}   
\label{intro}

\subsection{\bf Semidefinite programs and shortcomings of the usual  dual} 
Semidefinite Programs (SDPs) -- optimization problems with linear objective, linear constraints, and semidefiniteness constraints on matrix variables -- 
are some of the most versatile and popular optimization problems to emerge in the last thirty years. SDPs appear  
in combinatorial optimization, polynomial 
optimization, engineering, and other application areas,  and can be solved by efficient optimization algorithms. See, for example, 
\cite{NestNemirov:94, Ren:01} for the foundational theory of interior point methods, 
\cite{TohTodTut:99, sturm1999using} for efficient implementations of such methods, and \cite{yang2015sdpnal+, burer2003nonlinear, helmberg2000spectral, monteiro2024low} 
for efficient algorithms based on different principles.

We formulate an SDP mathematically as 

\begin{equation}\label{p}
	\begin{array}{rl} 
		\inf  & \,\, \la C, X   \ra  \\
		s.t. & \,\, \la A_i, X \ra \, = \, b_i \, (i=1, \dots, m) \\
		& \,\, X \succeq 0, 
	\end{array} \tag{$P$} 
\end{equation}
where the $A_i$ and $C$ are $n \times n$ symmetric matrices and $b \in \rad{m}. \, $ 
Also, for symmetric matrices 
$T$ and $S$ we write $S \preceq T$ to say that $T - S$ is positive semidefinite (psd) and 
$\la T, S \ra :=  {\rm trace}(TS)$ to denote their inner product. 

The problem \eqref{p}, which we call the {\em primal},  has a natural dual problem 
\begin{equation}\label{d}
	\begin{array}{rl} 
		\sup & \,\, \la b, y \ra  \\
		s.t.   & \,\,  \sum_{i=1}^m  y_i A_i \preceq C, 
	\end{array} \tag{$D$}
\end{equation}
where we write $\la b, y \ra$ for the inner product $b^\top y$  of $b$ and $y. \, $ In the examples 
it will be convenient to state the dual in terms of $C - \sum_{i=1}^m y_i A_i $ being psd. 
We will call this matrix a {\em slack matrix}. 

One of the most important roles of \eqref{d} is to certify boundedness of the optimal value of \eqref{p}    and optimality of feasible solutions.
For example, when $X$ is feasible in \eqref{p}, and $y$ in \eqref{d}, then the {\em weak duality}  inequality 
$\la C, X \ra \geq \la b, y \ra $ always holds. 
Thus, if we find a pair $X$ and $y$ whose objective values are equal, then we know they must be both optimal.

While weak duality is useful, we usually want a stronger property to hold, both for 
theoretical and for practical reasons. 
A desirable property of \eqref{p} and of \eqref{d} is {\em strong duality,} which is said to hold 
when the  optimal values of \eqref{p} and \eqref{d} agree, and the latter is attained, when it is finite. 
However, strong duality between \eqref{p} and \eqref{d} can fail, 
as the following example shows: 

\begin{Example} \label{example-main} 
	\begin{equation}    \label{problem-main}    
		\ba{rrcl} 
		\inf  &  \biggl \langle \underbrace{\bpx 0 & 1 & 0 \\ 
			1 & 0 & 0 \\ 
			0 & 0 & 0 \epx}_{C}, X \biggr \rangle   \\ \vspace{-.3cm} \\ 
		s.t.   & \biggl \langle \underbrace{\bpx 1 & 0 & 0 \\ 
			0 & 0 & 0 \\ 
			0 & 0 & 0 \epx}_{A_1}, X \biggl \rangle  & = & 0 \\ \vspace{-.3cm}  \\
		& \biggl \langle  \underbrace{\bpx 0 & 0 & 1 \\ 
			0 & 1 & 0 \\ 
			1 & 0 & 0 \epx}_{A_2}, X \biggl \rangle  & = &  0    \\ \vspace{-.3cm}  \\
		&  \biggl \langle \underbrace{\bpx 0 & 0 & 0 \\ 
			0 & 0 & 0 \\ 
			0 & 0 & 1  \epx}_{A_3}, X \biggl \rangle  & = &  1  \\ \vspace{-.3cm}  \\
		& X & \succeq & 0. 
		\ena 
	\end{equation} 
	
		We claim that its optimal value is zero. Indeed, assume  that  $X$ is feasible in it. Then by the first equation the $(1,1)$ element of $X$ is zero, and since $X$ is psd, 
	the first row and column of $X$ is zero. This proves our claim.
	
	The dual is 
	\begin{equation} \label{problem-main-dual} 
		\begin{aligned}
			\sup\quad & y_3 \\
			\text{s.t.}\quad &
			C - \sum_{i=1}^3 y_i A_i \, = \, \begin{pmatrix}
				-y_1 & 1    & -y_2 \\
				1   & -y_2 & 0    \\
				-y_2 & 0    & -y_3
			\end{pmatrix} \succeq 0. \\
		\end{aligned}
	\end{equation}		
	We claim that  \eqref{problem-main-dual} has no solution with value zero. Indeed, suppose $y$ is such a solution, so $y_3=0.$ Let us call the slack matrix $S.$ 
	Since $S \succeq 0, \, $ we deduce $y_2 = 0, \, $ a contradiction 
	to the $(1,2)$ element of $S$ being $1.$ 
	
	With some more work one can show that the optimal value of  \eqref{problem-main-dual}   is zero, but it is not attained.
	
\end{Example} 

Many examples of pathological SDPs are known, and most textbooks and surveys give such examples.   
Example \ref{example-main} is inspired by Theorem 3.1 in \cite{bomze2012think}: they show how to cleverly position a nonzero in $C,$ where 
the variable matrix $X$ is forced to be zero, in order  to create an instance with unattained dual optimal value. 

Strong duality can be ensured if we assume certain regularity conditions. 
The best known such   condition   is strict feasibility: 
when \eqref{p} is strictly feasible, i.e., it has a positive definite feasible solution,  then strong duality holds between \eqref{p} and \eqref{d}.
An analogous result holds when 
\eqref{d} is strictly feasible, i.e., when there is $y$ such that the slack matrix 
$C - \sum_i y_i A_i$ is positive definite.

However, assuming strict feasibility is not satisfactory from a theoretical perspective. Most importantly,  
it is of no help in finding an exact alternative system  of \eqref{p}, i.e., a semidefinite system which is feasible, exactly when 
\eqref{p} is infeasible. Indeed,  the  usual ``Farkas lemma"  system 
\begin{equation} 
	\begin{split}
		\sum_{i=1}^m y_i A_i  & \succeq  0, \, \,  
		\la b, y \ra  = -1 
	\end{split}\tag{$\mathit{alt \mhyphen P}$}  \label{altp} 
\end{equation}  
of \eqref{p} is not an exact alternative system   \footnote{More precisely, 
	\eqref{altp} {\em is}  an exact alternative system of \eqref{p}, if there is $y$ such that $\sum_i y_i A_i$ is positive definite.  
	However, this assumption is quite restrictive. We can of course also assume that all the $A_i$ and $C$ are diagonal, so \eqref{p} is just a linear program, but this assumption is even more restrictive.}: there are  instances of \eqref{p} which are infeasible, while \eqref{altp} is also 
infeasible. For a concise treatment of duality in conic linear programs, which include SDPs, see, e.g. \cite[Chapter 3]{Ren:01}, or 
\cite[Chapter 2]{BentalNem:01}.

\subsection{\bf Ramana's dual} 
Thirty years ago, in a seminal paper Ramana \cite{Ramana:97} \footnote{The first version of \cite{Ramana:97} was circulated in 1995.} 
constructed an elegant  dual problem, which avoids  the shortcomings of 
the traditional  dual.  Ramana's dual has the following striking  properties: i) it is a polynomial size explicit SDP ii) it assumes that \eqref{p} is feasible, but does not assume it is strictly feasible  iii) strong duality holds between \eqref{p} and Ramana's dual.
Put simply, it has all desirable properties of \eqref{d} when \eqref{p} is strictly feasible, without actually assuming that \eqref{p} 
is strictly feasible!

Ramana's dual  yields an exact alternative system of \eqref{p}, and fundamental results in complexity theory.
The most important of these are: 
\begin{enumerate}
		\item In the real number  model of computing, deciding feasibility of SDPs is in NP $\cap$ co-NP.
	\item In the Turing model of computing, deciding feasibility of SDPs is 
	\begin{enumerate}
		\item either in NP $\cap$ co-NP or not in  NP $\cup$ co-NP
		\item not NP-complete, unless 	NP = co-NP.
	\end{enumerate} 
\end{enumerate}
These results are still the best known on SDP feasibility. 

To state Ramana's dual, we assume that the primal \eqref{p} is feasible, and we denote by 
$\val()$ the optimal value of an optimization problem. 
We denote by $\sym{n}$ the set of $n \times n$ symmetric matrices, 
and by $\psdn$  the set of symmetric psd matrices.
We also introduce the linear operator $\A$ and its adjoint $\A^*$ as  
\begin{equation} \nonumber 
	\A X  = ( \la A_1, X \ra,  \dots,  \la A_m, X \ra)^\top,  \; 	\A^* y := \sum_{i=1}^m y_i A_i  \,  \text{for} \; X \in \symn, y  \in \rad{m}. 
\end{equation}

\begin{Theorem}\label{theo:main}
	Consider the optimization problem called the {\em Ramana dual of \eqref{p}:} 
	\begin{equation}  \label{d-ramana}  
		\begin{array}{rrcl}
			\sup        & \la b, y \ra  &     \\
			\text{s.t.} & C - \mathcal{A}^* y & \in & \mathcal{S}_+^n + \tan\left( U_{n-1} \right)  \\
			& y   & \in & \rad{m} \\
			& U_0 & = &  V_0 \, = \, 0                                     \\
			&  \multicolumn{3}{l}{
				\left.
				\begin{array}{rcl}
					\hspace{.5cm} \mathcal{A}^* y^i & \,  = \,    & U_i + V_i \\
					\la b, y^i \ra    & \, =  \,   & 0 \\
						y^i               & \, \in \, & \rad{m} \\
					U_i               & \, \in \,  & \mathcal{S}_+^n \\
					V_i               & \, \in \,  & \tan\left( U_{i-1} \right) 
								\end{array}
				\right\} \quad \text{for } i =1,\dots, n-1 
			} \\
		\end{array} \tag{D$_{\text{Ram}}$}
	\end{equation}
	Here for $U \in \psd{n}$ the set $\tan(U)$ is defined as 
	\begin{equation} \label{eqn-tanU-lambda} 
		\tan(U) = \biggl\{ W + W^\top \, : \, W \in \mymatrix{n \times n}, \, \bpx U & W \\
		W^\top & R \epx \in \psd{2n}  \; \text{for some} \, R \in \psd{n} \,   \biggr\}.  
	\end{equation}
	We then have 
	$$\val \eqref{p} = \val \eqref{d-ramana}, \,$$ and $\val \eqref{d-ramana}$ is attained when finite. 
	\qed 
\end{Theorem}
Note that in \eqref{d-ramana} the $y^1, \dots, y^{n-1}$ are vectors in $\rad{m}.$ To avoid confusion, 
we write $y_i$ for the $i$th component of the variable vector $y \in \rad{m}.$ 
     
We consider the $y$ the ''main" variable in \eqref{d-ramana}, since it plays a role analogous to the 
role of the $y$ variable in \eqref{d}. 
Thus we will usually say that $y$ is feasible in \eqref{d-ramana} with some $\{y^i, U_i, V_i \} $ 
and understand 
that the index $i$ runs from $0$ to $n-1. \, $ 	Also, in the examples we will not exhibit $U_0, V_0, \, $ and $V_1, \, $ since these are always zero:
$U_0 = V_0 = 0$ by definition, and $V_1 \in \tan(0) = \{0\}.$ 

\begin{Example} \label{example-main-cont}       
	(Example \ref{example-main} continued) 	In the Ramana dual of \eqref{problem-main} we claim  that 
	$y=0$ is a feasible  solution with value zero with some $\{y^i, U_i, V_i \}.$
	Indeed, let 
	$y^1 = e^1, \, y^2 = e^2, \,$ where here and in what follows,
	$e^i$ denotes the $i$th unit vector of appropriate dimension. 
	
	To construct the $U_i$ and $V_i$ we first observe that for $0 \leq r \leq n$ and 
	\begin{equation} \label{eqn-U-form} 
		U = \bpx I_r & 0 \\ 0 & 0 \epx \in \psd{n}
	\end{equation}
	any matrix in $\symn$ in which all nonzeros are in the first $r$ rows and columns is in $\tan(U).$
	Thus writing 
	\begin{equation} \label{eqn-UiVi-for-example} 
		\A^* y^1 = A_1 =  \underbrace{\bpx 1 & 0  & 0 \\ 0 & 0 & 0 \\  0 & 0 & 0 \epx}_{U_1}, \, \A^* y^2 = A_2 = \underbrace{\bpx 1 & 0  & 0 \\ 0 & 1 & 0 \\  0 & 0 & 0 \epx}_{U_2} + \underbrace{\bpx -1  & 0  & 1 \\ 0 & 0 & 0 \\  1 & 0 & 0 \epx}_{V_2}
	\end{equation}
	we see that $V_2 \in \tan(U_1)$ 
	and $C - \A^* y = C \in \tan(U_2) \subseteq \psd{3} + \tan(U_2), \, $ as wanted.
	
\end{Example}
In the following remarks we clarify the properties of \eqref{d-ramana}.
First,  the feasible set of \eqref{d-ramana} with respect to the $y$ variable is at least as large as the feasible set of \eqref{d}: this is because $0$ is in the tangent space of any psd matrix. Thus 
\begin{equation} \label{eqn-dram-vs-d} 
	\val \eqref{d-ramana} \geq \val \eqref{d}. 
\end{equation}
Second, in the next two useful formulas we connect the variables in \eqref{d-ramana} with the variables in \eqref{p}.
For that, suppose $y$ with some $\{y^i, U_i,V_i \}$ is feasible in \eqref{d-ramana}
and $X$ is feasible in \eqref{p}.  Then, as is standard in duality theory, we deduce 
\begin{equation} \label{eqn-CX-by} 
	\la C, X \ra - \la b, y \ra = \la C, X \ra - \la \A X, y \ra = \la C - \A^* y, X \ra. 
\end{equation}
Also, for $i=1, \dots, n-1$ we see that  
\begin{equation} \label{XUiVi} 
	\la X, U_i + V_i \ra = \la  X, \A^* y^i \ra = 	\la \A X, y^i \ra  =  \la b, y^i \ra  = 0.  
\end{equation}
Third, continuing the preceding argument, assume $X$ is actually a  strictly feasible solution in 
\eqref{p}. Since $V_1=0, $ from \eqref{XUiVi} we deduce 
\begin{equation} \label{eqn-XU1} 
	\begin{array}{rclrclrcl} 
		\la X, U_1 \ra = 0 & \Rightarrow & 	U_1 = 0 & \Rightarrow & V_2 = 0 & \Rightarrow & \la X, V_2 \ra = 0.
	\end{array}
\end{equation}

Repeating this argument with $U_2, \dots, U_{n-1}$ in place of $U_1, \, $ 
(when we start with $U_{n-1}$ we only need the first implication) 
we deduce 
that all $U_i$ and $V_i$ are zero, so in this case \eqref{d-ramana} is equivalent to 
\eqref{d} 
\footnote{This argument was inspired by a comment of Javier Pe{\~n}a, 
	whose help is gratefully acknowledged.}. 

Fourth, suppose now that \eqref{d-ramana} is strictly feasible. Then in particular, it has a feasible solution 
in which $U_1$ is positive definite. Hence by \eqref{eqn-XU1} we deduce the somewhat surprising conclusion that the only feasible solution of \eqref{p} is $X=0.$ 

Fifth, Ramana derived his dual for an inequality constrained SDP, i.e., 
for our dual \eqref{d}. 
A Ramana type dual for an equality constrained SDP, i.e., for our 
\eqref{p} was stated in \cite{RaTuWo:97}, and a result analogous to our Theorem \ref{theo:main} was also proved there. Our \eqref{d-ramana} has  some important differences with respect to 
the one stated in \cite{RaTuWo:97}, as follows. First, we isolated the $y$ variable to make clear 
that \eqref{eqn-dram-vs-d} holds. Second, we described the definition of 
the ``$\tan$" constraints separately in \eqref{eqn-tanU-lambda}, rather than plugging it into 
\eqref{d-ramana} as in several previous works. Third, we permit 
the $V_i$ to be in $\tan(U_{i-1}), \, $ whereas previous works restricted 
the $V_i$ to be in a subset of $\tan(U_{i-1}).  \, $
We focus on \eqref{d-ramana} in this particular form, since this form lends itself to a simple and intuitive analysis. 

Lastly, geometrically,  $\tan(U)$ is the tangent space of $\psd{n}$ at $U,$
defined as  
\begin{equation} \label{eqn-tandef-geo} 
	\tan(U) \, = \,  \Bigl\{ \, V \in \symn \, :  \,  \dist( U \pm \epsilon V, \psdn)   \rightarrow  \, \text{as} \, \epsilon \searrow   0   \, \Bigr\},
\end{equation}	
where $\dist(X, \psd{n}) \, = \, \inf  \{ \norm{X - Y} \, | \, Y \in \psd{n} \, \}$ is the distance 
of matrix \mbox{$X \in \symn$} from $\psd{n}.$ \co{For example, suppose  $U$ is as given in \eqref{eqn-U-form}, and all nonzeros of 
	$V \in \symn$ are in the first $r$ rows and columns. Then a simple calculation shows that 
	a perturbation of order $\epsilon^2$ 
	brings $U \pm \epsilon V$ back to $\psdn.$  }
However in what follows, we will rely only on the algebraic description of 
the tangent space given in \eqref{eqn-tanU-lambda}.

\subsection{\bf Literature} 	 Ramana's dual  is fundamental, however, the original proof of its correctness is  somewhat lengthy and  
technical. Thus  several papers gave shorter proofs, and explored connections to other work.  
Ramana, Tun\c{c}el and Wolkowicz \cite{RaTuWo:97} and \cite{Pataki:00B, Pataki:13, lourencco2023simplified} 
connected Ramana's dual to  the facial reduction algorithm of Borwein and Wolkowicz  \cite{BorWolk:81}. Klep and Schweighofer \cite{KlepSchw:12} 
designed a dual with similar properties, based on algebraic geometry.  Luo, Sturm, and Zhang \cite{LuoSturmZhang:97} gave a different proof of the correctness 
of Ramana's dual.  Ramana and Freund \cite{RaFreund:96} showed its usual Lagrange dual has the same optimal value as the original SDP. Generalizations are also available: 
\cite[Corollary 1]{Pataki:13} described a Ramana type dual for conic linear programs, 
assuming the underlying cone belongs to the class of {\em nice} cones. Further, \cite[Theorem 2]{liu2017exact} described a Ramana type dual for an arbitrary conic linear program.
These latter results are more general, however, they are also stated in a more abstract setting, so they 
have not led to complexity results comparable to Ramana's. 

Ramana's dual was used by de Klerk et al \cite{de2000self}  in 	self-dual embeddings.    
Due to its complexity implications it  is often mentioned in the discrete mathematics and theoretical computer science literature,  see for example, Lov{\'a}sz~\cite{lovasz2019graphs} and O' Donnell~\cite{o2017sos}.  Ramana's dual 
is often cited in surveys and books:  see for example, 	de Klerk \cite{de2006aspects}, Todd \cite{Todd:00}, 
Vandenberghe and Boyd  \cite{vandenberghe1996semidefinite}, Nemirovski \cite{nemirovski2007advances}, and    
Laurent and Rendl 	\cite{laurent2005semidefinite}. 

For completeness we list some references, which are less closely related, 
but also aim at understanding the complexities of SDP. 
Ramana's dual 	inspired many papers whose aim is to understand SDP duality, and the pathological phenomena 
that occur in it. The author in \cite{Pataki:17, pataki2019characterizing}   characterized badly behaved semidefinite {\em systems}, in which strong duality fails for some objective function.  Louren{\c{c}}o, Muramatsu, and Tsuchiya  in 
\cite{lourencco2021solving} showed how with a suitable oracle one can classify feasibility statuses of SDPs. 

Another stream of research addressed the issue of feasible solutions in SDPs, whose size (bitlength) 
is exponential in the size of the input.  The first such concrete example was constructed by Porkolab and Khachiyan \cite{porkolab1997complexity}. More recently, 
O' Donnell \cite{o2017sos} showed that such large solutions arise in sum-of-squares (SOS) 
proofs of nonnegativity. 
Raghavendra and Weitz \cite{raghavendra2017bit} and Gribling, Polak, and Slot \cite{gribling2023note} 
followed this line of research, and gave conditions that guarantee polynomial size solutions in  SOS proofs. 
Further, the author and Touzov \cite{pataki2024exponential} showed that large size solutions are more frequent than previously thought: 
 they arise in  SDPs with large so called singularity degree, after a simple linear transformation.

Despite the importance of Ramana's dual and the many followup papers,  one can make the case that we still need to understand it better. On  the one hand, 
the cited references characterize its optimal value. However, it would also be very useful to characterize its feasible set, 
both from the theoretical, and possibly a practical perspective.
Second, a simple correctness proof, accessible to both the optimization and the theoretical computer science communities, 
is also desirable.

\subsection{\bf Contributions}

We first connect Ramana's dual to a seemingly very different way of inducing strong duality: reformulating 
\eqref{p}  into a {\em rank revealing (RR) form} \cite{LiuPataki:15}, which helps us
verify the maximum rank of a feasible solution.
The RR form is constructed using elementary row operations (inherited from Gaussian elimination),
and rotations. 
Second, while previous works characterized its optimal value, 
here we completely characterize its feasible set. In particular 
we show it is a higher dimensional representation, or {\em lift}  of 
a dual problem with similar favorable properties, which, however is not an explicit SDP. 
Thus, our work provides a connection to the theory of lifts, 
representations of 
optimization problems in a higher dimensional space: see for example the recent 
survey \cite{fawzi2022lifting}. 
\co{As a corollary, we obtain a short and transparent derivation of Ramana's dual, which we believe  is accessible to both the optimization 
	and the theoretical computer science communities.

	Second, in Theorem \ref{thm:connect:feasibleset}   we completely characterize the feasible set of \eqref{d-ramana}. Third, using elementary linear algebra, we similarly derive 
	the Ramana dual of \eqref{d}. 
}

As a corollary, 
we obtain a short and elementary  proof of Theorem \ref{theo:main}, and of its counterpart Theorem \ref{theo:main:dual}, which derives the Ramana dual of 
\eqref{d}; we hope our proofs will be accessible to both the optimization and theoretical computer science communities.

As we mentioned in the abstract, our approach is combinatorial. While 
a \mbox{``combinatorial approach"}   is not perfectly defined, the main features of our proofs are:
\begin{enumerate}
	\item \label{1111} We avoid the use of  most concepts in convex analysis, such as relative  interiors, faces, and conjugate faces, which play 
	an important role in the analysis of \cite{RaTuWo:97, Pataki:00B, Pataki:13}. 
	In fact, we only use a single ingredient from continuous optimization theory, a theorem of the alternative, which we state as a proposition for convenience:
	\begin{Proposition} \label{proposition-alternative} 
		Suppose \eqref{p} is feasible. Then it is not strictly  feasible $\Leftrightarrow$ 
		the system 
		\begin{equation} \label{eqn-A*y-alternative} 
			\A^* y \in \psd{n} \setminus \{ 0 \}, \,  \la b, y \ra = 0
		\end{equation}
		is feasible 	\footnote{In turn, this result 	can be proved by the standard strong duality result between \eqref{p} and \eqref{d}, assuming strict feasibility in one of them.}.
	\end{Proposition}

	\item \label{2222} We show that feasible solutions in \eqref{d-ramana} are identified with regular facial reduction sequences, i.e., essentially discrete structures. 
\end{enumerate}

\subsection{\bf Organization of the paper and guide to the reader} In Subsection 
\ref{subsection-preliminaries} we fix notation, prove three  simple propositions, and define one of the main players of the paper,
regular facial reduction sequences.
In 	Section  \ref{section-d-ramana}  we analyse \eqref{d-ramana}: 
\begin{itemize}
	\item  In Subsection \ref{subsection-certify-slack} we recall the rank revealing (RR)  form of \eqref{p} from \cite{LiuPataki:15}. 
	This  form makes it easy to verify 	the maximum rank of a feasible matrix in \eqref{p}. We then  show how to construct the RR form. 
	\item In Subsection \ref{subsection-strongdual} we study  the  strong dual  of \eqref{p},  which has all the properties required from \eqref{d-ramana}. However, 
	the strong dual  relies on knowing a maximum rank feasible solution in \eqref{p},
	and such a solution in general is not known known explicitly. 
	\item In Subsection \ref{subsection-connect} in Theorem
	\ref{theo:liftdstrong} 
	we give our first characterization  the feasible set 
	of \eqref{d-ramana}: we show it is a higher dimensional representation, a lift,
	of the feasible set of 
	the strong dual. As a corollary, we prove  Theorem \ref{theo:main}. 
	\item Ramana's dual may look somewhat magical at first, so in Subsection 
	\ref{subsect-intuition} we give intuition how it naturally arises from the RR form and the strong dual. 
	\item 	While in Subsection \ref{subsection-connect} we described the "$y$" portion of 
	feasible solutions of \eqref{d-ramana}, this is not yet a complete characterization, as it does not characterize the $\{y^i, U_i, V_i \}$ portion of feasible solutions. 
To complement Subsection \ref{subsection-connect}, in Subsection \ref{subsection-characterize-feasible-set} 
	we completely characterize its feasible set. We believe that such a characterization
	is essential for a potential succesful implementation.
\end{itemize}
While the results of Subsections \ref{subsection-certify-slack} and \ref{subsection-strongdual}  are known, the proofs in this paper are much 
simpler, and, as we alluded before,
rely on much  less machinery than the proofs in \cite{LiuPataki:15}. Theorem 
\ref{thm:connect:feasibleset} is related to Corollary 1 in \cite{Pataki:13}. In that result we considered a conic linear program stated a so-called {\em nice cone}, and characterized the dual cone of the so-called {\em minimal cone} of 
\eqref{d}. That result, however, is stated in a more technical manner, whereas in the main part of the current paper we do not refer  to dual cones, or minimal cones.
All the other results are new. 

We complete the paper with Appendix \ref{section-ramana-primal}, where 
  we derive corresponding results for the Ramana dual of \eqref{d}. 
These  results follow from results from Section \ref{section-d-ramana} and some elementary linear algebra, so most of them are only sketched.  
\co{In Appendix 
\ref{section-connect-FR} we explain the connection to classical versions of facial reduction algorithms. 
}

We organized the paper's results to be accessible to a broad audience. Some readers  may  only want to see a quick and transparent derivation of 
\eqref{d-ramana}. For them,  reading only Subsection \ref{subsection-preliminaries}, and Section \ref{section-d-ramana}, until, and including the proof of Theorem 
\ref{theo:main} will  suffice. 

\subsection{\bf Preliminaries}    
\label{subsection-preliminaries} 

We denote by $\sym{n,k}$ the set of $n \times n$ symmetric matrices in which all nonzeroes 
appear in the first $k$ rows and columns.
We let $\psd{n,k} = \psdn \cap \sym{n,k}$ i.e, the set of psd matrices in which only the 
upper left $k \times k$ block can be nonzero. 
We denote by 
$\pd{n,k}$ the matrices in $\psd{n,k}$ in which the upper left $k \times k$ block is positive definite.

Pictorially, $U$ in equation \eqref{eqn-tan-shape} is in $\psd{n,k}$ and 
$V$ is in $\sym{n,k}.$  In this equation and later 
$\oplus$ stands for a psd submatrix, and the $\ti$ stands for a block with arbitrary elements.

\vspace{.1cm}
\begin{equation} \label{eqn-tan-shape}
	U \, = \, \begin{pmatrix}
		\bovermat{$k$}{{\mbox{$\,\, \oplus \,\,$}}}	& \bovermat{$ n-k $}{{\mbox{$\,\,\,0\,\,\,\,\,$}}}  	\\ 
		0 & 0 
	\end{pmatrix}, V \, = \, \begin{pmatrix}
		\bovermat{$k$}{{\mbox{$\,\, \ti \,\,$}}}	& \bovermat{$ n-k $}{{\mbox{$\,\,\,\ti \,\,\,\,\,$}}}  	\\ 
		\ti & 0 
	\end{pmatrix}.
\end{equation}

Next we state three   basic propositions. 
The proofs of Proposition \ref{prop-XYQ} and \ref{prop-UVQ} are
straightforward from the properties of the trace and 
the definition of $\tan(U).$ 

\begin{Proposition} \label{prop-XYQ} 
	Suppose $Q$ is an $n \times n$ orthonormal matrix. Then 
	\begin{equation}  \label{eqn-XY} 
		\ba{rcl}
		\la  S, T  \ra & = & \la Q^\top S Q, Q^\top T Q \ra 
		\ena 
	\end{equation}
	for all $S, T \in \symn.$   
	Further, 
	\begin{equation} \label{eqn-UV} 
		V \in \tan(U) \; \Leftrightarrow \; Q^\top V Q \in \tan(Q^\top U Q)
	\end{equation}
	
	for all $U \in \psdn, \, $ and $V \in \symn.$ 
	\qed 
\end{Proposition}

\begin{Proposition}  \label{prop-UVQ} The following hold: 
	\begin{enumerate}
		\item 	 \label{prop-UVQ-psd} If $U \in \psd{n,k}$ and $V \in \tan(U), \, $ then  $V \in \sym{n,k}. \, $ 
		\item 	 \label{prop-UVQ-pd}  If $U \in \pd{n,k}$ and $V \in \sym{n,k}, \, $ then $V \in \tan(U). \, $  
	\end{enumerate}
	\qed 
\end{Proposition}
We can visualize Proposition \ref{prop-UVQ} in equation 
\eqref{eqn-tan-shape}.
If $U$ is as given on the left, and $V \in \tan(U), \, $ then $V$ must be of the  form given on the right. Further, if the $\oplus$ block in $U$ is positive definite, then 
any $V$ in the form on the right is in $\tan(U).$

\begin{Proposition} \label{prop-maxrank-in-C} 
	Suppose $C$ is a convex subset of $\symn$ and $X$ is a maximum rank psd matrix in $C$ of the form 
	\begin{equation} \label{eqn-Y-form}  
		X =  
		\bpx 0 & 0 \\
		0 & \Lambda \epx, 
	\end{equation}
	where $\Lambda$ is order $r$ and positive definite. Then in any psd matrix in $C$ the first $n-r$ rows and columns are zero.
\end{Proposition}	
\proof{} Let us denote the nullspace of any  matrix $B$ by $\N(B).$ 
Assume to the contrary that $X' $ is a psd matrix in $C$ 
and the first $n-r$ rows and columns of $X'$ are not all zero.
Let $X'' = \frac{1}{2} (X + X').$ We then claim
$$
\N(X'') = \N(X) \cap \N(X') \subsetneq \N(X).  
$$
Indeed, the equality is from basic 
linear algebra \footnote{If $X \in \psdn, \, u \in \rad{n},$ then 
	$u \in \N(X) \Leftrightarrow u^\top X u = 0.$}. Also, the $\subsetneq$ relation holds,  since $\N(X) \setminus \N(X')$ is nonempty  
(for example any vector whose last $r$ elements are zero, but is not in 
$\N(X')$  is in this set).
Thus, $X''  \in C$ and   has larger rank than $X, \, $ a contradiction. \qed

The following notation will be useful. If $r_0, \dots, r_t$ are real numbers,
$0 \leq k \leq \ell \leq t, \, $ then we write 
$$
r_{k:\ell} := \sum_{i=k}^\ell r_i.
$$
For brevity, we omit parantheses in this notation: for example, we write $r_{1:i+1}$ instead of 
$r_{1:(i+1)}$

We next introduce a main player of the paper: 
\begin{Definition} \label{definition-regfr} 
	We say that $Y_1, \dots, Y_k$ is a {\em regular facial reduction sequence}\footnote{Slightly different versions of regular facial reduction sequences have been defined in other papers, e.g. in \cite{pataki2022echelon}.}  {\em for $\psd{n}$} 
	if the  $Y_i$ are  in $\symn$ and are of the form

		$$
		Y_1   = 
		\bordermatrix{
			& \overbrace{\qquad}^{\textstyle r_{1}} &  \overbrace{\qquad \qquad}^{\textstyle{n-r_1 }} \cr\\
			& \Lambda_1  &  0   \cr
			& 0   &  0  \cr}, \, \dots, Y_i  = 
		\bordermatrix{
			& \overbrace{\qquad }^{\textstyle r_{1:i-1}} & \overbrace{\qquad}^{\textstyle r_{i}} & \overbrace{\qquad}^{\textstyle n-r_{1:i}} \cr\\
			& \times  &  \times  &  \times \cr
			& \times  &  \Lambda_i    &  0 \cr
			& \times  &  0  &  0 \cr} 
		$$
	for  $i= 1,\dots, k. \,$ Here the $r_i$ are nonnegative integers, the $\Lambda_i$ diagonal positive definite matrices, and the $\times$ symbols 
	correspond to blocks with arbitrary elements.
	
\end{Definition}

Note that regular facial reduction sequences are essentially discrete structures. 
When we use them, we only use  that the $\Lambda_i$ are positive definite, and what their sizes are;
however, we never refer to their actual entries.

\section{Analysis of \eqref{d-ramana}} 
\label{section-d-ramana} 

\subsection{\bf The rank revealing (RR) form of \eqref{p} and reformulations}  \label{subsection-certify-slack} 

\begin{Definition} \label{definition-RR-P} 
	We say that \eqref{p} is in rank revealing form, or RR form, if for some $0 \leq k \leq m$ 
	\begin{enumerate}
		\item \label{definition-RR-D-1} 
		$A_1, \dots, A_k$ is a regular facial reduction sequence in which the sizes of the positive definite blocks are 
		nonnegative integers $r_1, \dots, r_k, \, $ respectively. 
		\item \label{definition-RR-D-2}
		$b_1 = \dots = b_k = 0.$ 
		\item  \label{definition-RR-D-3}  there is a feasible solution of the 
		form
		\begin{equation} \label{eqn-Y-form} 
			\bpx 0 & 0 \\
			0 & \Lambda \epx, 
		\end{equation}
		in \eqref{p}, where $\Lambda$ is order 	$n - r_{1:k}, \,$ and positive definite. 
	\end{enumerate} 
	If \eqref{p} is in RR form, then we also say that the first $k$ equations in \eqref{p} certify  that the solution 
	in \eqref{eqn-Y-form} has maximum rank in it. For brevity, sometimes we say that 
	the first $k$ equations in \eqref{p} certify the maximum rank.
\end{Definition}
We next explain the  terminologies in Definition  \ref{definition-RR-P}.  
Suppose  \eqref{p} is in RR form, as given in Definition \ref{definition-RR-P}, and the sizes of the positive definite blocks in $A_1, \dots, A_k$ are $r_1, \dots, r_k, \, $ respectively.
Also suppose $X$ is feasible in \eqref{p}.  Since the upper left order $r_1$ block of $A_1$ is positive definite, by $\la A_1, X \ra = 0$ we deduce  the corresponding 
block of $X$ is zero. Since $X$ is psd, the first $r_1$ rows and columns of $X$ are zero. Then  $\la A_2, X \ra = 0$ implies 
the next  $r_2$ rows and columns of $X$ are zero; etc. Thus the first $k$ equations indeed certify that the solution given in \eqref{eqn-Y-form}  has maximum rank 
in \eqref{p}. 

\begin{Example} \label{example-main-cont2} 
	(Example \ref{example-main} continued) We claim the SDP \eqref{problem-main} is in RR form without any reformulation. 
	Indeed, $(A_1, A_2)$ is a regular facial reduction sequence, and the first two equations certify that the unique feasible solution 
	$$
	X = 
	\bpx 0 & 0 & 0 \\ 
	0 & 0 & 0 \\ 
	0 & 0 & 1  \epx
	$$
	has maximum rank.
	
\end{Example}

Next we look at how to transform \eqref{p} into RR  form, if it is not in that form to start with.

\begin{Definition} \label{definition-reformulation} 
	We say that we 
	\begin{enumerate}
		\item {\em rotate} a set of matrices say $M_1, \dots, M_k$ by an orthonormal matrix $Q, \,$ if we replace 
		$M_i$ by $Q^\top M_i Q$ for all $i.$ We  say that we rotate \eqref{p} by an orthonormal matrix $Q$ if we rotate all $A_i$ and 
		$C$ by $Q.$ 
		
		\item 
		{\em reformulate}  \eqref{p}  if we apply the following operations (in any order):
		\begin{enumerate}
			\item  We rotate all $A_i$ and $C$ by an orthonormal matrix. 
			\item \label{exch}  For some $i \neq j$ we exchange equations 
			$$\la A_i, X \ra = b_i \; \text{and} \; \la A_j, X \ra = b_j. \;  $$
			\item \label{replace} We replace an equation by a linear combination of equations. That is, for some $i \in \{1, \dots, m \}$ 
			we replace 
			$$\la A_i, X \ra = b_i \; \text{by} \; \la \A^* y, X \ra = \la b, y \ra \; \text{where} \; y \in \rad{m} \; \text{ and } \; y_i \neq 0.  $$
		\end{enumerate}
		
		\item We say that  by reformulating \eqref{p} we obtain a reformulation. 
		
	\end{enumerate}
\end{Definition}

Note that operations \eqref{exch} and \eqref{replace} in Definition \ref{definition-reformulation}
are elementary row operations inherited from Gaussian elimination.

As the next lemma shows, the simple operations of Definition \ref{definition-reformulation} suffice to put \eqref{p} into RR form.
\begin{Lemma} \label{lemma-RR-P} 
The SDP \eqref{p} can always be reformulated into RR form.
\end{Lemma}

\proof{} If \eqref{p} is strictly feasible, then we do not have to reformulate it, we just set $k=0.$ 
If \eqref{p} is not strictly feasible, then we invoke Proposition \ref{proposition-alternative} and find  $y \in \rad{m}$  
such that 
\begin{equation} \label{eqn-A*y-alternative-2} 
\A^* y \in \psd{n} \setminus \{ 0 \}, \,  \la b, y \ra = 0. 
\end{equation}
Let $Q$ be a matrix of orthonormal eigenvectors of $\A^* y$ and assume w.l.o.g. that the first element of
$y$ is nonzero. Replace
$(A_1 , b_1)$ by  
$(\A^* y, 0), $ then rotate all $A_i$ by $Q.$ After this we have 
$$
A_1 = \bpx \Lambda_1 & 0 \\ 0 & 0 \epx,
$$
where $\Lambda_1$ is diagonal and positive definite, of order, say $r_1.$

Next,  from \eqref{p} we construct a new SDP, say $(P')$ by deleting the first $r_1$ rows and 
columns from all $A_i$ and from $C.$ We see that \eqref{p} is equivalent to $(P'),$ since   
in any $X$ feasible of \eqref{p} the first  $r_1$ rows and columns must be zero. 
Thus we proceed in like fashion with $(P').$ 
\qed

Note that the construction in Lemma \ref{lemma-RR-P} is theoretical.
While the proof is constructive, to actually compute the RR form we would need to find 
$y$ feasible in \eqref{eqn-A*y-alternative-2}, and for that, we would need to solve an SDP 
in exact arithmetic.

\begin{Example} \label{example-bad-dual} 
Consider an SDP with data 
\begin{equation}  \label{problem-bad-n=4} 
	\ba{rclrclrcl} 
	A_1 & = & \left(\begin{array}{rrrr}  
		-4 & 15 & 6 & 3 \\
		15 & 3 & 0 & 5 \\
		6 & 0 & 5 & 0 \\
		3 & 5 & 0 & 0
	\end{array} \right),
	& 
	A_2 & = &  \left(\begin{array}{rrrr} 
		-1 & 6 & 2 & 1 \\
		6 & 1 & 0 & 2 \\
		2 & 0 & 2 & 0 \\
		1 & 2 & 0 & 0
	\end{array} \right), & 
	A_3 & = & \left(\begin{array}{rrrr}  
		2 & 3 & 0 & 0 \\
		3 & 0 & 0 & 1 \\
		0 & 0 & 1 & 0 \\ 
		0 & 1 & 0 & 0
	\end{array} \right), \\ \\ 
	C & = & \left(\begin{array}{rrrr}  
		1 & 0 & 0 & 0 \\
		0 & 1 & 0 & 0 \\
		0 & 0 & 1 & 0 \\ 
		0 & 0 & 0 & 0 
\end{array} \right), &    
b & =  &  (5, 2, 1)^\top.
\ena 
\end{equation}

Suppose we reformulate this SDP by performing the operations
\begin{equation} \label{eqn-for_bad-into-RR} 
	\begin{array}{rcl}
		(A_1, b_1) & = & (A_1, b_1) - 3 (A_2,b_2) + (A_3,b_3), \\
		(A_2, b_2) & = & (A_2,b_2) - 2 (A_3,b_3). 
	\end{array}
\end{equation}
We thus obtain the SDP with data 
\begin{equation}  \label{problem-good-n=4} 
	\ba{rclrclrcl} 
	A_1 & = & \left(\begin{array}{rrrr}  
		1 & 0 & 0 & 0 \\
		0 & 0 & 0 & 0 \\
		0 & 0 & 0 & 0 \\
		0 & 0 & 0 & 0
	\end{array} \right),
	& 
	A_2 & = &  \left(\begin{array}{rrrr} 
		-5 & 0 & 2 & 1 \\
		0 & 1 & 0 & 0 \\
		2 & 0 & 0 & 0 \\
		1 & 0 & 0 & 0
	\end{array} \right), & 
	A_3 & = & \left(\begin{array}{rrrr}  
		2 & 3 & 0 & 0 \\
		3 & 0 & 0 & 1 \\
		0 & 0 & 1 & 0 \\ 
		0 & 1 & 0 & 0
	\end{array} \right), \\ \\ 
	C & = & \left(\begin{array}{rrrr}  
		1 & 0 & 0 & 0 \\
		0 & 1 & 0 & 0 \\
		0 & 0 & 1 & 0 \\ 
		0 & 0 & 0 & 0
	\end{array} \right), 
	&    
	b & =  &  (0, 0, 1)^\top.
	\ena 
\end{equation}
We claim that this SDP  is in RR form. Indeed,
\begin{equation}\label{eqn-X-maxrank}
	X = \left(\begin{array}{rrrr}
		0 & 0 & 0 & 0 \\
		0 & 0 & 0 & 0 \\
		0 & 0 & 1 & 0 \\
		0 & 0 & 0 & 1
	\end{array}\right)
\end{equation}
is feasible in it. Also, the matrices $(A_1, A_2)$ form a regular facial reduction sequence
(with $r_1 = r_2 = 1$), which certify that $X$ has maximum rank. (In fact, 
$(A_1, A_2,A_3)$ is also a regular facial reduction sequence, but $A_3$ does not play a role in certifying the maximum rank, since $b_3 \neq 0.$)

Naturally, the $X$ above  is also the maximum rank feasible solution 
in the system defined by the original $A_i$ and $b$ in \eqref{problem-bad-n=4}. However,  from this form of the $A_i$ and $b$  this would be difficult to tell.

\end{Example} 

Lemma \ref{lemma-RR-P} describes a  kind of {\em facial reduction algorithm:} in an RR form 
the 
first $k$ constraints force any feasible $X$ to have its first $r_{1:k}$ rows and columns equal to zero, 
i.e., to live in a face of $\psd{n}$ \footnote{A convex subset $F$ of $\psd{n}$ is a face of $\psd{n},$ 
	if $X, Y \in  F$ and $\frac{1}{2}(X+Y) \in F$ imply 
	that $X$ and $Y$ are both in $F.$ The faces of $\psdn$ are exactly the sets $T^\top \psd{n,k} T, \,$ for some $k \in \{0, \dots, n \}, $ and an invertible matrix $T$ 
	\cite{Pataki:00A} }. Facial reduction algorithms originated in \cite{BorWolk:81}, then simpler versions were introduced by 
Waki and Muramatsu \cite{WakiMura:12} and the author \cite{Pataki:00B, Pataki:13}. 
Our treatment in this paper is sufficiently simplified that we do not even have to define faces.

It is clear that we can always construct an RR form with $k \leq n,$ 
since we can drop any equation $\la A_i, X \ra = 0$ in which 
$r_i, \,$ the size of the positive definite block is zero.
The next lemma shows that we can do a bit better.

\begin{Lemma} \label{lemma-kleqn-1} 
There is always an RR form with $k \leq n-1.$ 
\end{Lemma}
\proof{} Suppose \eqref{p} is in RR form as given in Definition \ref{definition-RR-P}.
Then $k \leq n$ follows, as we argued above.

Suppose $k = n.$ We claim that in this case there is an RR form with $k=1.$ 
Indeed, if $\lambda_1  > 0$ is sufficiently large then in  $A_2' := \lambda_1 A_1 + A_2$ the upper left order $2$  block is positive definite:
this follows by the Schur-complement condition for positive definiteness. 
Similarly, if  $\lambda_2  > 0$ is sufficiently large then in $A_3' := \lambda_2 A_2' + A_3$ the upper left order $3$  block is positive definite. Continuing, 
we construct an equation $\la A_n', X \ra = 0$ with $A_n'$ positive definite, so after a rotation we indeed 
obtain  an RR form with $k=1$ (and the only feasible solution being $X=0$). 
\qed

Next we discuss how reformulating \eqref{p} affects feasible solutions of 
\eqref{p}, \eqref{d} and \eqref{d-ramana}. For that, we note that a reformulation of \eqref{p} can be encoded just by two matrices, say $M$ and $Q$ as follows.
The elementary row operations amount to replacing 
$\A$ by $M \A$ and $b$ by $M b, \, $ where $M \in \mymatrix{m \times m}$ is invertible.
Also, to construct the  reformulation we can just use $Q, \, $ the 
product of all rotation matrices used in the reformulation process.

The proof of the following proposition is straightforward from \eqref{eqn-UV} in 
 Proposition \ref{prop-XYQ}.

\begin{Proposition} \label{prop-reform-invariance} 
	Suppose we reformulate \eqref{p} and the reformulation is represented by 
	matrices $M$ and $Q$ as described above. Then 
	\begin{enumerate}
		\item \label{prop-reform-invariance-p}  $X$ is feasible in \eqref{p} before the reformulation iff 
		$Q^\top X Q$  is feasible after the reformulation.
		\item \label{prop-reform-invariance-d} $y$ is feasible in \eqref{d} before the reformulation iff 
		$M^{-*}y$   is feasible  after the reformulation. 
		\item \label{prop-reform-invariance-d-ramana} $y$ with $\{y^i, U_i, V_i \}$ is feasible in \eqref{d-ramana} 
		before the reformulation iff 
		$M^{-*}y$  with $\{M^{- *} y^i,$ $ Q^\top U_i Q,$ $ Q^\top V_i Q \}$   is feasible  after the reformulation
		\footnote{Here, and in what follows, for a linear operator $M$ we write $M^{- *}$ for the inverse of the adjoint $M^*.$ }.
	\end{enumerate}
	\qed 
\end{Proposition}

\subsection{\bf The strong dual of \eqref{p}}
\label{subsection-strongdual} 
In this subsection we first state a strong dual of \eqref{p}, which has the same number of variables as \eqref{d}, but has all the properties we require from Ramana's dual. 

\begin{Lemma}  \label{lemma-strong-dual} 
Suppose a maximum rank solution 
in \eqref{p} is of the form
\begin{equation} \label{eqn-max-rank-form}  
	Q \bpx 0 & 0 \\
	0 & \Lambda \epx Q^\top, 
\end{equation}
where $Q$ is orthonormal, and $\Lambda$ is order $r$ and positive definite.

Consider the optimization problem called the {\em strong dual  of \eqref{p},} 
\begin{equation}\label{d-strong}
	\begin{split}
		\sup  & \,\, \la b, y \ra    \\
		s.t. & \,\,  C - \A^* y = Q V Q^\top   \\ 
		& \;\;\;\;\;\;\;\;\;\;\;\;V \in \symn, V_{22} \in \psd{r} 
	\end{split} \tag{\mbox{${\rm D_{{\rm strong},Q}}$}}
\end{equation}
where $V_{22}$ stands for the lower right order $r$ block of $V.$ 
We then have 
\begin{equation} \label{eqn-val-D-val-re-P} 
	\val \eqref{p} = \val \eqref{d-strong}, 
\end{equation} 
and $\val \eqref{d-strong}$ is attained when finite.
\qed
\end{Lemma}


\proof{} Suppose we rotate \eqref{p} by $Q.$ By Proposition \ref{prop-reform-invariance}  
we see that $X$ is feasible before the rotation iff 
$Q^\top X Q$ is feasible afterwards. Thus the optimal value and attainment in \eqref{p} does not change by this rotation. Also, 
$y \in \rad{m}$ is feasible in \eqref{d-strong} before the rotation 
iff it is feasible in (\mbox{${\rm D_{{\rm strong},I}}$})
after  the rotation. Thus we will assume without loss of generality that $Q=I.$ 


Let $(P')$ be the SDP obtained from \eqref{p} by deleting the first $n-r$ rows and columns 
in all $A_i$ and in $C, \, $  and 
$(D')$ the dual of $(P').$ 
We   claim that 
\begin{equation} \label{eqn-all-equal} 
\val \eqref{p} = \val (P') = \val (D') = \val \eqref{d-strong},
\end{equation} 
and that the optimal values of $(D')$ and \eqref{d-strong} are attained. 

Indeed,  in the first equality $\leq$ follows, since by Proposition \ref{prop-maxrank-in-C} 
in any $X$  feasible solution of \eqref{p} the first $n-r$ rows and columns are zero.
In the same inequality $\geq$ follows, since if $X'$ is feasible in $(P')$ then adding $n-r$  all zero rows and columns gives a feasible solution of \eqref{p}.
The second equality in \eqref{eqn-all-equal} follows, since $(P')$ is strictly feasible. Strict feasibility in $(P')$ also implies 
that $\val (D')$ is attained. The last equation and attainment in \eqref{d-strong} follow, since the feasible set of $(D')$ and 
\eqref{d-strong} are the same. 

\qed

\vspace{.5cm} 
Note that the slack matrix $C - \A^* y$ 
in feasible solutions of \eqref{d-strong} is of the form 

\begin{equation} \label{eqn-S-d-strong} 
	C - \A^* y =  Q \begin{pmatrix}
		\bovermat{$n-r$}{{\mbox{$\,\, \,\,\ti \,\,\,\,$}}}	& \bovermat{$ r  $}{{\mbox{$\,\,\,\,\ti \,\,\,\,\,$}}}  	\\ 
		\ti & \oplus
	\end{pmatrix} Q^\top ,
\end{equation} 
where, as usual, the blocks marked by $\ti$ contain arbitrary elements and the $\oplus$ block is positive semidefinite.

\begin{Example} (Example \ref{example-main}  continued) 
	Consider again the SDP \eqref{problem-main} in which the maximum rank feasible solution is 
	$$
	X = 
	\bpx 0 & 0 & 0 \\ 
	0 & 0 & 0 \\ 
	0 & 0 & 1  \epx.
	$$
	Thus in the strong dual we can take $Q=I.$ We repeat the usual dual 
	here from \eqref{problem-main-dual} for convenience: 
		\begin{equation} \label{problem-main-dual-repeat} 
		\begin{aligned}
			\sup\quad & y_3 \\
			\text{s.t.}\quad &
			C - \A^* y \, = \, \begin{pmatrix}
				-y_1 & 1    & -y_2 \\
				1   & -y_2 & 0    \\
				-y_2 & 0    & -y_3
			\end{pmatrix} \succeq 0. \\
		\end{aligned}
	\end{equation}	
		The  strong dual is 
	just like the usual dual \eqref{problem-main-dual-repeat}, 
	except only the lower right $1 \times 1$ corner of the slack matrix 
	must be psd, i.e., nonnegative. Hence  
	$y=0$ is feasible (and optimal) in this strong dual. 
	
	Thus, as expected from  Lemma \ref{lemma-strong-dual}, 
	strong duality holds between \eqref{problem-main} and its strong dual.

\end{Example}

\begin{Example} \label{example-bad-dual-continued}  (Example \ref{example-bad-dual}  continued) 
	Consider the SDP with data \eqref{problem-good-n=4}. 
	We saw that the first two rows and columns 
	of any feasible $X$ are zero, hence the third constraint implies $x_{33}=1. \, $ 
	Thus the optimal value is $1.$ 
	
	The dual is
	\begin{equation} \label{problem-example-bad-dual-continued-dual} 
		\begin{aligned}
			\sup \quad & y_3 \\
			\text{s.t.}\quad &
			C - \A^* y \; = \;
			\begin{pmatrix}
				1 - y_1 + 5y_2 - 2y_3 & -3y_3 & -2y_2 & -y_2 \\
				-3y_3 & 1 - y_2 & 0 & -y_3 \\
				-2y_2 & 0 & 1 - y_3 & 0 \\
				- y_2 & - y_3 & 0 & 0
			\end{pmatrix} \succeq 0, \\
			& y_1,y_2,y_3 \in \mathbb{R},
		\end{aligned}
	\end{equation}
	and it is clear that in any feasible solution $y_3 = 0. \,$ Hence 
	there is a positive duality gap.
	
	Let us next examine  the strong dual. For that, we observe that the SDP has a maximum rank solution \eqref{eqn-X-maxrank} in which the lower right $2 \times 2$ block is positive definite, and the other elements are zero. Thus in the strong dual we can take $Q=I, \,$ 
	and in the slack matrix $C - \A^* y$ only the lower right $2 \times 2$ block
	must be psd. 
	Thus any $y$ with $y_3=1$ is feasible, and optimal in the strong dual. 
	
	 Hence, again, as expected from  Lemma \ref{lemma-strong-dual}, 
	strong duality holds between the SDP and its strong dual.

\end{Example}

Example \ref{example-bad-dual-continued} also shows that the RR form can help us verify when 
strong duality fails between \eqref{p} and \eqref{d}. Indeed,  we saw that 
the gap between the optimal values of the SDP defined by \eqref{problem-good-n=4}, which is in RR form,
and its dual is $1. \, $ Thus  by 
Proposition \ref{prop-reform-invariance} the same is true of the SDP defined by 
\eqref{problem-bad-n=4}, which is {\em not} in RR form, and its dual. 
However, this latter statement would be much more difficult to verify directly.

\co{Given that \eqref{d-strong} already achieves strong duality, and it has just  
the same number of variables as the original dual 
\eqref{d}, we may wonder, why we even need \eqref{d-ramana}?
Note however, that \eqref{d-strong} relies on knowing the maximum rank feasible solution in 
\eqref{p}, which, as we discussed before, is not readily available: to compute it, we would need to solve SDPs in exact arithmetic.
}

\subsection{\bf The feasible set of \eqref{d-ramana} as a lift of the feasible set of 
	\eqref{d-strong}} 
\label{subsection-connect}

Next we state one of the main results of the paper. It shows that we can use 
the strong dual as a building block to construct \eqref{d-ramana}.
\begin{Theorem} \label{theo:liftdstrong} 
	There is a $Q$ orthonormal matrix with the following properties:
	\begin{enumerate}
		\item \label{theo:liftdstrong-1} A maximum rank solution in \eqref{p} is of the form given in \eqref{eqn-max-rank-form}.
		\item \label{theo:liftdstrong-2} For any $y \in \rad{m}$ it holds that 
		\begin{equation*}
			y \text{ is feasible in } \eqref{d-strong}
			\Leftrightarrow 
			y \text{ is feasible in } \eqref{d-ramana} \text{ with some } \{y^i, U_i, V_i\}.
		\end{equation*}

	\end{enumerate}
	\qed
\end{Theorem}  
Before we get to the proof, we explain Theorem \ref{theo:liftdstrong}.
Its essence is that the feasible set of 
\eqref{d-strong} is the projection of the feasible set of \eqref{d-ramana}, which lives in a 
higher dimensional space. 

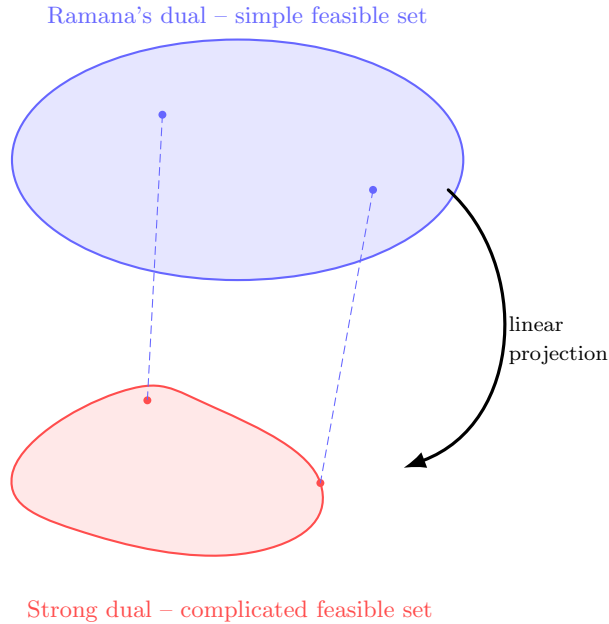
\begin{figure}[h]
	\centering
	\begin{tikzpicture}[line join=round, line cap=round, >=Latex, scale=1]
		
		\fill[blue!10] (0,4) ellipse (3 and 1.6);
		\draw[thick, blue!60] (0,4) ellipse (3 and 1.6);
		
		\fill[blue!60] (-1.0,4.6) circle (1.5pt);
		\fill[blue!60] ( 1.8,3.6) circle (1.5pt);
		
		\path
		(-3,-0.2)  coordinate (L1)
		(-1.6,0.9) coordinate (L2)
		(-0.3,0.7) coordinate (L3)
		( 1.1,-0.3) coordinate (L4)
		( 0.3,-1.2) coordinate (L5)
		(-1.9,-1.0) coordinate (L6);
		
		\fill[red!9]
		plot[smooth cycle, tension=0.8] coordinates { (L1) (L2) (L3) (L4) (L5) (L6) };
		\draw[thick, red!70]
		plot[smooth cycle, tension=0.8] coordinates { (L1) (L2) (L3) (L4) (L5) (L6) };
		
		\fill[red!70] (-1.2,0.8) circle (1.5pt) coordinate (pA);
		\fill[red!70] ( 1.1,-0.3) circle (1.5pt) coordinate (pB);
		
		\draw[blue!60, densely dashed] (-1.0,4.6) -- (pA);
		\draw[blue!60, densely dashed] ( 1.8,3.6) -- (pB);
		
		\draw[very thick, -{Latex[length=3mm]}]
		(2.8,3.6) .. controls (3.9,2.6) and (3.8,0.5) .. (2.2,-0.1)
		node[midway, right, align=left, inner sep=2pt] {\footnotesize linear\\[-1pt]\footnotesize projection};
		
		\node[blue!60, align=center] at (0,5.9)
		{\small Ramana's dual -- simple feasible set};
		
		\node[red!70, align=center] at (-0.1,-2.0)
		{\small Strong dual -- complicated feasible set};
		
	\end{tikzpicture}
	\caption{The feasible set of Ramana's dual projects onto the feasible set of the strong dual (stylized).}
	\label{fig:projection} 
\end{figure}

Why do we even need this higher dimensional representation? 
While \eqref{d-strong} already achieves strong duality, its feasible set has a 
complicated description, as it needs to know a maximum rank feasible solution in 
\eqref{p}. Such a maximum rank feasible solution in 
\eqref{p},  as we discussed before, is not readily available: to compute it, we would need to solve SDPs in exact arithmetic.
On the other hand, \eqref{d-ramana} has many more variables, but it has a favorable representation, as it is an explicit SDP.
See Figure \ref{fig:projection} (inspired by a figure in \cite{fawzi2022lifting})) for a schematic view. 

Thus, by recent terminology, the feasible set of \eqref{d-ramana} is a {\em lift} of the feasible set of \eqref{d-strong}: we refer to  \cite{fawzi2022lifting} for a survey of lifts -- a beautiful area in convex optimization. Lifts of polyhedra also appear in combinatorial optimization, and go by the name of 
{\em extended formulations:} for a recent survey we refer to 
\cite{conforti2010extended}.

\proof{We first prove \eqref{theo:liftdstrong-1}.
	Let $X$ be a maximum rank solution in \eqref{p}, and let $r$ denote its rank. 
	Consider a reformulation of \eqref{p} into RR form, and let $Q$ be the product of all rotation matrices in the reformulation process. Then $Q^\top X Q$ is a maximum rank solution 
	after the reformulation, in which the lower right order $r$ block is positive definite, and all other elements are zero. Thus \eqref{theo:liftdstrong-1} holds. 
	
	To prove \eqref{theo:liftdstrong-2} suppose we rotate \eqref{p} by $Q.$ 
	As we discussed in the proof of Lemma \ref{lemma-strong-dual},
	$y \in \rad{m}$ is feasible in \eqref{d-strong} before the rotation 
	iff it is feasible in (\mbox{${\rm D_{{\rm strong},I}}$})
	after  the rotation. Also, the second statement in \eqref{theo:liftdstrong-2}
	is invariant under this rotation by \eqref{prop-reform-invariance-d-ramana} in Proposition \ref{prop-reform-invariance}.  
	Thus without loss of generality 
	we assume $Q=I.$ 
	
	We start with the implication 
	$\Rightarrow, $ so suppose $y$ is feasible in \eqref{d-strong}.
	Suppose that in an RR form of \eqref{p} the first $k$ equations certify the maximum rank solution.
	By Lemma \ref{lemma-kleqn-1} we assume $k \leq n-1.$ 
	Since these $k$ equations are a linear combination of the original equations, there is 
	$y^1, \dots, y^k \in \rad{m}$ such that 
	\begin{eqnarray} \label{eqn-A*y1A*yn-1} 
		(\A^* y^1, \dots, \A^* y^k) 
		&       & \text{is a regular facial reduction sequence, and}  \\ \label{eqn-byi}
		\langle b, y^i \rangle = 0 
		&  & \text{for } i = 1, \dots, k. 
	\end{eqnarray}
	For $i=1, \dots, k$ let  $\Lambda_i$ be  the positive definite block in 
	$\A^* y^i $ and let 
	$r_i$ be the order of $\Lambda_i.$ Then we decompose $\A^* y^i$ into $U_i + V_i$ as 
	
	\begin{equation} \label{eqn-define-Ui-gabor} 
		\underbrace{	\begin{pmatrix}
				\bovermat{$r_{1:i-1}$}{\,\,\,\,\,\ti \,\,\,\,\,} & \bovermat{$\,\,\, r_{i}$}{\,\,\,\, \ti \,\,\,\,}	& \bovermat{$n - r_{1:i}$}{\,\,\,\,\,\,\, \ti \,\,\,\,\,\,\,}	\\ 
				\ti  & \, \Lambda_{i}   \,  & 0  \\ 
				\ti    &0  & 0 
		\end{pmatrix}}_{\A^* y^i} =  
		\underbrace{	\begin{pmatrix}
				\bovermat{$r_{1:i-1}$}{\,\,\,\,\,I \,\,\,\,\,} & \bovermat{$\,\,\, r_{i}$}{\,\,\,\, 0\,\,\,\,}	& \bovermat{$n - r_{1:i}$}{\,\,\,\,\,\,\, 0 \,\,\,\,\,\,\,}	\\ 
				0  & \, \Lambda_{i}  \,  & 0  \\ 
				0    &0  & 0 
		\end{pmatrix}}_{U_i}  + \underbrace{\begin{pmatrix}
				\bovermat{$r_{1:i-1}$}{\,\,\,\,\,\ti \,\,\,\,\,} & \bovermat{$\,\,\, r_{i}$}{\,\,\,\, \ti \,\,\,\,}	& \bovermat{$n - r_{1:i}$}{\,\,\,\,\,\,\, \ti \,\,\,\,\,\,\,}	\\ 
				\ti  & \, 0  \,  & 0  \\ 
				\ti    &0  & 0 
		\end{pmatrix}}_{V_i},
	\end{equation}
	i.e., we simply define $V_i := \A^* y^i - U_i.$ 
	In \eqref{eqn-define-Ui-gabor}  we indicated the blocks with arbitrary elements by $\ti$ marks. 
	Thus we have 
	\begin{equation} \label{eqn-A*yidecomp} 
		\left. \begin{array}{rcl}
			\hspace{.5cm} \A^* y^i   & \,  = \,    & U_i + V_i \\
			U_i               & \, \in \,  & \mathcal{S}_+^n \\
			V_i               & \, \in \,  & \tan\left( U_{i-1} \right)
		\end{array}
		\right\} \quad \text{for } i =1,\dots, k. 
	\end{equation} 
	Since $k \leq n-1$ we need to ``pad" the sequence $\{y^i, U_i, V_i \}$ 
	with zeros. That is, we add $n-1-k$ to the index of each, and  
	define $y^i = 0$ and 
	$U_i = V_i = 0$ for $i=1, \dots, n-1-k.$ Then  \eqref{eqn-byi}  and \eqref{eqn-A*yidecomp}  
	hold with $n-1$ in place of $k.$
	
	To complete the proof, we see that $r_{1:n-1} = n -r, \, $ hence    
	$U_{n-1} \in \pd{n,n-r}. \, $ Since $y$ is feasible in the  strong dual of \eqref{p}, the lower right 
	order $r$ block of  $C - \A^* y $  is psd.   (Recall that now $Q=I.$)
		Thus, 
	\begin{equation} \label{eqn-CA*ytan} 
		C - \A^* y \in \psd{n} + \tan(U_{n-1}).
	\end{equation}
	Combining all of the above 
	completes the proof. 
	
	To prove the $\Leftarrow$ implication in \eqref{theo:liftdstrong-2}, 
	suppose $y$ with some $\{y^i, U_i, V_i\}$ is feasible in \eqref{d-ramana}. 
	Recall from \eqref{XUiVi} that $\la X, U_i + V_i \ra = 0$ for $i=1, \dots, n-1.$ 
	Thus, repeating the argument in \eqref{eqn-XU1}  almost verbatim, we get 
	\begin{equation} \label{eqn-XU1XV2} 
		\begin{array}{rclrclrcl} 
			\la X, U_1 \ra = 0 & \Rightarrow & 	U_1 \in \psd{n, n - r} & \Rightarrow & V_2 \in \sym{n, n-r}  & \Rightarrow & \la X, V_2 \ra = 0,
		\end{array}
	\end{equation}
	where the first implication follows since the lower right order $r$ corner of $X$ is positive definite, so this block of $U_1$ is zero, and by 
	$U_1 \succeq 0.$ The second implication is by $V_2 \in \tan(U_1)$ and by part \eqref{prop-UVQ-psd} of Proposition \ref{prop-UVQ}.
	Repeating this with $U_2, \dots, U_{n-1}$ in place of $U_1, \, $ 
	(when we do it with  $U_{n-1}$ we only need the very first implication) 
	we see that 
	all $U_i$ are in $\psd{n, n - r},$ so  
	$$ \tan(U_{n-1}) \subseteq \sym{n, n - r}.$$ 
	Hence the lower right order $r$ block of $C - \A^*y$ is psd. 
	Thus $y$ is feasible in \eqref{d-strong}, as wanted. 
}

\qed 

Now  we can prove Theorem \ref{theo:main}. Let $Q$ be as in Theorem 
\ref{theo:liftdstrong}. Then 
\begin{equation} \nonumber 
	\val \eqref{p} = \val \eqref{d-strong} = \val \eqref{d-ramana},
\end{equation} 
where the first equation is from Lemma \ref{lemma-strong-dual} and the second is from 
Theorem \ref{theo:liftdstrong}. Also by Lemma \ref{lemma-strong-dual}, the optimal value of 
\eqref{d-strong} is attained, hence by Theorem \ref{theo:liftdstrong} the optimal value of  \eqref{d-ramana} is also attained. Thus the proof is complete.
\qed

\begin{Example} (Example \ref{example-bad-dual} continued)
	Consider again the SDP with data  \eqref{problem-good-n=4}. 
	We will construct an optimal solution for  its Ramana dual. For that, first let us fix 
an arbitrary $y \in \rad{3}$ whose last element is $1.$ Recall from Example \ref{example-bad-dual-continued}
that 	$y$ is optimal in the strong dual.

	We then need 
	a suitable  
	$\{y^i, U_i,V_i\}.$ To construct the $y^i$ (which are in $\rad{3}$),
	we note that this SDP is in RR form, so according to the proof 
	of Theorem \ref{theo:liftdstrong} we can take 
		\begin{equation}
		\begin{array}{rcl} 
			y^1 & = & 0, \\
			y^2 & = & e^1,  \\
			y^3 & = & e^2.
		\end{array} 
	\end{equation}
	We then  decompose the $\A^* y^i$ as 
	\begin{equation} \label{eqn-A*y2etc}
	\begin{array}{rcl}  
		\A^* y^1 & = & \underbrace{0}_{U_1 \succeq 0}   \\
		\A^* y^2 & = & \underbrace{\left(\begin{array}{rrrr}  
				1 & 0 & 0 & 0 \\
				0 & 0 & 0 & 0 \\
				0 & 0 & 0 & 0 \\
				0 & 0 & 0 & 0
			\end{array} \right)}_{U_2 \succeq 0} + \underbrace{0}_{V_2 \in \tan(U_1)} \\[1ex]
		\A^* y^3 & = & \left(\begin{array}{rrrr} 
			-5 & 0 & 2 & 1 \\
			0 & 1 & 0 & 0 \\
			2 & 0 & 0 & 0 \\
			1 & 0 & 0 & 0
		\end{array} \right)
		= \underbrace{\left(\begin{array}{rrrr} 
				1 & 0 & 0 & 0 \\
				0 & 1 & 0 & 0 \\
				0 & 0 & 0 & 0 \\
				0 & 0 & 0 & 0
			\end{array} \right)}_{U_3}
		+ \underbrace{\left(\begin{array}{rrrr} 
				-6 & 0 & 2 & 1 \\
				0 & 0 & 0 & 0 \\
				2 & 0 & 0 & 0 \\
				1 & 0 & 0 & 0
			\end{array} \right)}_{V_3 \in \tan(U_2)}.
	\end{array}
\end{equation}	
Since the lower right $2 \times 2$ block of $C - \A^* y$ is psd (in fact zero), 
we see that 
it is in $\psd{4} + \tan(U_3).$
Thus, $y$ with the $\{y^i, U_i, V_i\}$ is indeed feasible (and optimal) in the Ramana dual.

Similarly, we can also construct an optimal solution to the Ramana dual of the SDP defined by the 
original constraints \eqref{problem-bad-n=4}. For that, we  note it is brought into RR form by 
the 	operations listed in	\eqref{eqn-for_bad-into-RR}, and no rotation. 
So again, let us take any $y$ whose last element is $1,$ and 
according to the proof 
of Theorem \ref{theo:liftdstrong}, let  
\begin{equation}
	\begin{array}{rcl} 
		y^1 & = & 0, \\
		y^2 & = & (1,-3,1)^\top, \\
		y^3 & = & (0,1,-2)^\top,
	\end{array} 
\end{equation}
and a decomposition listed in \eqref{eqn-A*y2etc}. 
\co{

	Similarly,  we can construct an optimal solution to the Ramana dual of the SDP given by data in \eqref{problem-good-n=4}. Given that this SDP is in RR form, this task is much simpler: now we can take $y^1 = 0, y^2 = e^1, y^3 = e^2,$ and the decomposition the same as in 
	\eqref{eqn-A*y2etc}. }We leave the details to the reader.
	
\end{Example}

\subsection{\bf Remarks for better intuition} 
\label{subsect-intuition}

Ramana originally derived his dual using very different arguments from ours, 
and he derived it for our dual \eqref{d}. 
The original result in his paper \cite{Ramana:97} as well as correctness of 
our \eqref{d-ramana} may look like ``magic" at first, so 
in this subsection we explain the intuition behind it. 

To derive \eqref{d-ramana}, we need two ingredients: the RR form, and the strong dual. 
\begin{enumerate}
	\item 

The  
RR form arises very naturally: it is just an iterated 
version  of the classical theorem of the alternative given in Proposition 
\ref{proposition-alternative}, combined with some basic linear algebra. 
\item 
The second ingredient, the strong dual \eqref{d-strong} is also  natural: once we know 
what the maximum rank solutions in \eqref{p} look like, 
its correctness follows since 
the restricted primal $(P')$ is strictly feasible. See the proof of Lemma \ref{lemma-strong-dual}.
\item 
Given these ingredients, we still need to create an explicit SDP from them.
For that, we first observe that the matrices in an RR form and the slack matrix 
$C - \A^* y$ in the strong dual naturally decompose into 
a psd part and a tangent space part.
Second, the tangent space of the semidefinite cone is 
representable by psd constraints. Third, the key relation 
$$
V_i \in \tan(U_{i-1}) 
$$ is preserved by rotations. 
\end{enumerate}

\subsection{\bf A complete characterization of the feasible set of \eqref{d-ramana}}
\label{subsection-characterize-feasible-set} 

The implication $\Leftarrow$ in part \eqref{theo:liftdstrong-2} of Theorem \ref{theo:liftdstrong} 
characterizes the $"y"$ part of feasible solutions of \eqref{d-ramana}. However, it is not a complete characterization, as it 
does 
not characterize the $\{y^i, U_i, V_i \}$ portion. In the next result we complement 
Theorem \ref{theo:liftdstrong} and completely characterize
the feasible solutions of \eqref{d-ramana}.

\begin{Theorem} \label{thm:connect:feasibleset}  
Suppose $y$ with some $\{y^i, U_i, V_i\}$ is feasible in \eqref{d-ramana}. 
Then after a suitable rotation of \eqref{p} the following holds:
\begin{enumerate}
	\item  \label{thm:connect:feasibleset-1}  $(\A^* y^1, \dots, \A^* y^{n-1})$ is a regular facial reduction sequence.
	\item \label{thm:connect:feasibleset-2}  If the size of the positive definite block in $\A^* y^i$ is $r_i$ for all $i, \,$ 
	then   $$U_1 \in \psd{n,r_1}, U_2 \in \psd{n,r_{1:2}}, \dots, U_{n-1} \in \psd{n,r_{1:n-1}}.$$ 
\end{enumerate}

\end{Theorem} 
\proof{} Let us make the assumption. For brevity, let $Y_i := \A^* y^i$ for all $i.$ 
We will  rotate all  $A_i, C, Y_i, U_i, V_i$  several times to achieve 
\eqref{thm:connect:feasibleset-1} and \eqref{thm:connect:feasibleset-2}. 
We first rotate all these matrices to achieve 
\begin{equation} \nonumber 
Y_1 = \bpx \Lambda_1 & 0 \\ 0 & 0 \epx, 
\end{equation}
where $\Lambda_1$ is diagonal positive definite. 
This can be done since $Y_1 = U_1 \in \psdn.$ Let $r_1$ be the order of $\Lambda_1.$ 

For the induction step, suppose that $1 \leq i < n-1$ and  the following invariant conditions  hold:
\begin{enumerate}[label=(inv-\arabic*), ref=inv-\arabic*, leftmargin=*]
\item \label{invariant-1} $Y_1, \dots, Y_i$ is a regular facial reduction sequence, in which the positive definite blocks have 
order $r_1, \dots, r_i,  \,$ respectively.
\item \label{invariant-2} $U_1 \in \psd{n,r_1}, U_2 \in \psd{n,r_{1:2}}, \dots, U_{i} \in \psd{n,r_{1:i}}.$ 
\end{enumerate}
Both these statements hold when $i=1.$ 
We will next make sure they hold with $i+1$ in place of $i.$ 
We have 
\begin{equation} \label{eqn-Yi+1Ui+1Vi+1}
Y_{i+1} = \underbrace{U_{i+1}}_{\in \psd{n}} + \underbrace{V_{i+1}}_{\in \tan(U_{i})}. 
\end{equation} 
By $U_{i} \in \psd{n,r_{1:i}}$ and  
\eqref{prop-UVQ-psd} in Proposition \ref{prop-UVQ} we deduce 
$V_{i+1} \in \sym{n,r_{1:i}}, \,$ so the lower right order 
$n -r_{1:i}$ 
block of $Y_{i+1}, \, $ which we call 
$\bY, \, $ is psd. 

We let $r_{i+1}$ be the rank of $\bY$ and $Q'$ be a matrix of orthonormal eigenvectors of $\bY \, $ 
and  deduce 

\vspace{.5cm} 

$$
\ba{rcl} 
Q^\top Y_{i+1} Q & = & 
\begin{pmatrix}
\bovermat{$r_{1:i}$}{\,\,\,\,\;\;\;\; \ti \,\,\,\,\,\;\;\;\;} & \bovermat{$r_{i+1}$}{\ti}	& \bovermat{$n - r_{1:i+1}$}{\,\,\,\,\,\,\,\,\,\,\,\,\, \times \,\,\,\,\,\,\,\,\,\,\,\,\,\,}	\\ 
\ti   &  \Lambda_{i+1}  & 0  \\ 
\times   &0  & 0 
\end{pmatrix} \!\!\!     
\end{array} \, , \; \text{where} \; Q = \bpx I_{r_{1:i}} & 0 \\
0	  & Q'  \\
\epx,
$$
and $\Lambda_{i+1}$ is diagonal, positive definite, and of order $r_{i+1}.$ 
So we rotate all matrices by $Q, \,$ and afterwards item \eqref{invariant-1} holds with $i+1$ in place of $i.$ Further,  \eqref{invariant-2} still holds with $i.$ Hence  we still have 
$V_{i+1} \in \sym{n,r_{1:i}}.$

Since all nonzeros in both $Y_{i+1}$  and $V_{i+1}$ are in the 
first $r_{1:i+1}$ rows and columns, by \eqref{eqn-Yi+1Ui+1Vi+1} 
 the same is true of 
$U_{i+1}.$ So \eqref{invariant-2} holds with $i+1$ in place of $i, $ as wanted. 
After we achieved \eqref{invariant-1} and \eqref{invariant-2} for $i=1, \dots, n-1,$ 
the proof is complete. 
\qed
   
\begin{Example} (Example \ref{example-bad-dual} continued)   
One way to illustrate Theorem \ref{thm:connect:feasibleset} is to again consider \eqref{p} defined by the data in \eqref{problem-bad-n=4}. 
For that, we recall from Example  \ref{example-bad-dual-continued} that 
we constructed an optimal solution to its Ramana dual.
\co{
, and saw 
that $(\A^* y^1, \A^* y^2, \A^* y^3)$ was a regular facial reduction sequence.
}
Now we construct another feasible solution to its Ramana dual, so we let 
\begin{eqnarray*} 
y & = & (1, 1, 0)^\top, \\
y^1 = y^2 & = & 0,  \\
y^3 & = & e^1.   
\end{eqnarray*}
We claim that $y$ with $y^1, y^2, y^3$ and some suitable $U_i$ and $V_i$ is a feasible 
solution. Indeed, this follows, since  
\begin{equation}
\ba{rrclrcl}
\A^* y^1 = \A^* y^2 \, = \, 0,  
& 
\A^* y^3  & = & \underbrace{\left(\begin{array}{rrrr}  
		\phantom{-}1 & \phantom{-}0 & \phantom{-}0 & \phantom{-}0 \\
		\phantom{-}0 & \phantom{-}0 & \phantom{-}0 & \phantom{-}0 \\
		\phantom{-}0 & \phantom{-}0 & \phantom{-}0 & \phantom{-}0 \\
		\phantom{-}0 & \phantom{-}0 & \phantom{-}0 & \phantom{-}0
	\end{array} \right)}_{U_3},
& 
C - \A^* y  & = &  \left(\begin{array}{rrrr}
	5 & 0 & -2 & -1 \\  
	0 & 0 & 0 & 0 \\  
	-2 & 0 & 1 & 0 \\  
	-1 & 0 & 0 & 0 
\end{array}\right)
\ena
\end{equation}
\co{$$ 
C - y_1 A_1 - y_2 A_2 - y_3 A_3 \;=\; C - A_1 - A_2 
= \begin{pmatrix}  
	5 & 0 & -2 & -1 \\  
	0 & 0 & 0 & 0 \\  
	-2 & 0 & 1 & 0 \\  
	-1 & 0 & 0 & 0  
\end{pmatrix}.
$$
}
thus $C - \A^* y \in \psd{4}+\tan(U_3).$ 

According to the proof of Theorem \ref{thm:connect:feasibleset}, after a suitable rotation $(\A^* y^1,$ $ \A^* y^2, \A^* y^3)$ becomes a regular facial reduction sequence and we can see that 
in this case no rotation is needed.
\end{Example} 

Theorem \ref{thm:connect:feasibleset} is of interest for two reasons. The first is theoretical: since Ramana's dual is fundamental, and 
many references characterized its  optimal {\em value},  it is also of interest to characterize its feasible set. 

The second is possibly practical: 
such a characterization is essential to successfully implement Ramana's dual.  Indeed, suppose 
a solver delivers a (possibly approximate)  solution to \eqref{d-ramana}. Then the proof of
Theorem \ref{thm:connect:feasibleset} shows how to construct the rotations to check whether the 
$\A^* y^i$ are n the right form. 

Note that a full scale implementation may be difficult (due to the large number of extra variables.). 
However, even a limited implementation, 
just using a few extra $y^i, U_i, \, $ and $V_i$ can fix the pathologies in several SDPs: this is true in theory, using exact arithmetic in all computations. 	It would be interesting to see whether this theory translates into a practical advance for SDP solvers. 

\subsection{\bf Ramana's exact alternative system} 

Ramana in \cite{Ramana:97} described an exact alternative system for \eqref{d} with the following two key features:
\begin{itemize}
	\item This system has the same data as the feasible set of \eqref{d}, namely the $A_i$ and $C.$
	\item  It is feasible exactly when \eqref{d} is infeasible.
\end{itemize}
In this subsection we describe an exact alternative system for \eqref{p}, in the spirit of Ramana's work.
Given that most proofs are straightforward modifications of proofs in the previous part of the paper, we only sketch most of them.

To motivate it, we first give an example: 
\begin{Example} \label{example-infeas} 
	Consider the semidefinite system \eqref{mot-ex}:
	\beq \label{mot-ex}
	\ba{cclcc}
	\underbrace{\bpx 1 & 0 & 0 \\
	0 & 0 & 0 \\
	0 & 0 & 0 
	\epx}_{A_1}  &\bullet&  X & = & \underbrace{0}_{b_1} \\
	\vspace{0.2cm}
	\underbrace{\bpx 0 & 0 & 1 \\
	0 & 1 & 0 \\
	1 & 0 & 0 
	\epx}_{A_2} \, &\bullet& \, X & = & \underbrace{-1}_{b_2}.  \\
	&& X & \succeq & 0,
	\ena
	\eeq
	We claim it is infeasible. Indeed, suppose $X$ is feasible in it, and let us write $x_{ij}$ for the 
	$(i,j)$ element of $X.$ By the first constraint we have $x_{11}=0$ and by psdness we see that the first row and column of 
	$X$ is zero. Thus the second constraint implies $x_{22}=-1, \, $ a contradiction.
	
	Yet, the traditional alternative system \eqref{altp} fails to certify infeasibility of 
	\eqref{mot-ex}: there is no $y \in \rad{2}$ such that $y_1 A_1 = y_2 A_2 \succeq 0, \, y_1 b_1 + y_2 b_2 = -1.$
	
\end{Example}

The main result of this subsection follows:
\begin{Theorem} \label{thm:Ram-alt} 
The SDP \eqref{p} is infeasible $\Leftrightarrow$ the semidefinite system 
\eqref{alt-ramana} below, called Ramana's alternative system is feasible: 
\begin{equation}  \label{alt-ramana}  
	\begin{array}{rrcl}
		        & \mathcal{A}^* y & \in & \mathcal{S}_+^n + \tan\left( U_{n-1} \right)  \\
		        & \la b, y \ra & = & - 1 \\
		& y   & \in & \rad{m} \\
		& U_0 & = &  V_0 \, = \, 0                                     \\
		&  \multicolumn{3}{l}{
			\left.
			\begin{array}{rcl}
				\hspace{.5cm} \mathcal{A}^* y^i & \,  = \,    & U_i + V_i \\
				\la b, y^i \ra    & \, =  \,   & 0 \\
				y^i               & \, \in \, & \rad{m} \\
				U_i               & \, \in \,  & \mathcal{S}_+^n \\
				V_i               & \, \in \,  & \tan\left( U_{i-1} \right) 
			\end{array}
			\right\} \quad \text{for } i =1,\dots, n-1 
		} \\
	\end{array} \tag{$\mathit{alt \mhyphen Ram \mhyphen P}$} 
\end{equation}

\end{Theorem}

\begin{Example} (Example \ref{example-infeas} continued)
	Let $y = e^2.$ We claim that this $y$ with a suitable $\{y^i, U_i, V_i \}$ is feasible 
	in the Ramana alternative system of \eqref{mot-ex}. Indeed, let 
	$y^1 = 0, \, y^2 = e^1,  \,$ and we show the decomposition of the $\A^* y^i$ and $\A^* y$ 
	below:  
\begin{equation} \label{eqn-UiVi-for-example} 
	\A^* y^1 = \underbrace{0}_{U_1 \succeq 0}, \, \A^* y^2 = \underbrace{\bpx 1 & 0  & 0 \\ 0 & 0 & 0 \\  0 & 0 & 0 \epx}_{U_2 \succeq 0} + \underbrace{0}_{V_2 \in \tan(U_1)}, \, \A^* y  = \underbrace{\bpx 1 & 0  & 0 \\ 0 & 1 & 0 \\  0 & 0 & 0 \epx}_{\succeq 0} + \underbrace{\bpx -1  & 0  & 1 \\ 0 & 0 & 0 \\  1 & 0 & 0 \epx}_{\in \tan(U_2)}
\end{equation}
\end{Example} 

One way to derive \eqref{alt-ramana} is by using Ramana's dual, \eqref{d-ramana}.
Here we give a derivation which we believe to be more concise, and elegant.

The first ingredient in our derivation is 
a theorem of the alternative analogous to the one stated in Proposition \ref{proposition-alternative}:

	\begin{Proposition} \label{proposition-alternative-not-strict} 
	The SDP \eqref{p} is not strictly  feasible $\Leftrightarrow$ 
	the system 
	\begin{equation} \label{eqn-A*y-alternative-notstrict} 
		\A^* y \in \psd{n} \setminus \{ 0 \}, \,  \la b, y \ra \leq  0
	\end{equation}
	is feasible.
	
	\qed 
\end{Proposition}
Note that Proposition \ref{proposition-alternative-not-strict} does not distinguish between SDPs which are feasible, just not strictly feasible;  and SDPs which are infeasible.
\co{\begin{Example}
	The semidefinite system 
	Take for example, the feasible system 
	\begin{equation}
		\begin{array}{rcl}
			\bpx 1 & 0 \\
			0 & 0 \epx \bullet X & = & 0 \\[2.5ex]  
			X & \succeq & 0,
		\end{array}
	\end{equation}
	is feasible. Suppose we create two new problem as follows: in the first, we add a constraint $x_{12}=1, \,$ which makes it infeasible; in the second, we add a constraint $x_{22}=1,$  which makes it feasible. 
	then 
\end{Example}
}

\begin{Example}
Proposition \ref{proposition-alternative-not-strict} produces the same certificate 
$$
y = (1,0)^\top
$$
for the feasible, but not strictly feasible system
\begin{equation}
	\begin{array} {rcl}
		\bpx 1  & 0 \\
		0  & 0 \epx \bullet X & = & 0 \\ [2.5ex]
		\bpx 0  & 0 \\ 
		0  & 1 \epx \bullet X & = & 1 \\
		X & \succeq & 0.
	\end{array} 
\end{equation}
and for the infeasible system 
\begin{equation}
	\begin{array}{rcl}
		\bpx 1 & 0 \\
		0 & 0 \epx \bullet X & = & 0 \\[2.5ex]  
		\bpx 0 & 1 \\
		1 & 0 \epx \bullet X & = & 1 \\
		X & \succeq & 0. 
	\end{array}
\end{equation}
\end{Example}

The following lemma is a counterpart of Lemma \ref{lemma-RR-P}:
\begin{Lemma} \label{lemma-infeas-certificate} 
	The SDP \eqref{p} is infeasible $\Leftrightarrow$ it has a reformulation in which for some 
	$0 \leq k \leq m$ the following hold: 
		\begin{enumerate}
			\item \label{prop-infeas-certificate-1}  
			$A_1, \dots, A_{k}$ is a regular facial reduction sequence.
			\item \label{prop-infeas-certificate-2}
			$b_1 = \dots = b_{k-1} = 0, \, b_k = -1.$ 
				\end{enumerate} 
		\end{Lemma}
{\em Proof sketch:} The easy direction $\Leftarrow$ is just like the discussion after the proof of 
 Lemma \ref{lemma-RR-P}. To get a contradiction, we assume $X$ is feasible in the reformulation.
 The first $k-1$ equations prove that the first $r_{1:k-1}$ rows and columns of $X$ are zero.
 Then the $k$th equation proves the trace of its diagonal block is $-1, \,$ the required contradiction.
 
 The more difficult direction $\Rightarrow$ follows by repeated application of Proposition 
 \ref{proposition-alternative-not-strict}, and noting that at some point we must find a $y$ 
 feasible in \eqref{eqn-A*y-alternative-notstrict} for which $\la b, y \ra < 0.$ 
 \qed
 
Now it is straightforward that whenever \eqref{p} is infeasible, there is such a reformulation with 
$k \leq n-1, \, $ since we can drop any one of the first $k$ equations in which the size of the positive definite block is zero.

{\em Proof sketch of Theorem \ref{thm:Ram-alt}:} 
$\Leftarrow:$ Suppose \eqref{alt-ramana} is feasible.
To get a contradiction, assume that 
\eqref{p} is also feasible, and let $X$ be a feasible solution in \eqref{p}. 
By rotating \eqref{p}, we assume that 
\begin{equation} \label{eqn-Y-form}  
	X =  
	\bpx 0 & 0 \\
	0 & \Lambda \epx, 
\end{equation}
for some $\Lambda$ positive definite matrix of order, say $r.$ 
 Using an argument like in Proposition 
\ref{prop-reform-invariance}, we see that after this rotation 
\eqref{alt-ramana} is still feasible. 
By the same argument as in \eqref{XUiVi}, we deduce that 
$$
\la X, U_i + V_i \ra = 0 \; \text{for} \, i=1, \dots, n-1. 
$$
We also repeat the argument in 
\eqref{eqn-XU1XV2} and deduce that all $U_i$ are in $\psd{n, n - r}, \, $ hence 
the lower $r \times r$ block of $\A^* y$ is psd. 
We then arrive at the contradiction 
$$
0 \leq \la X, \A^* y \ra \, = \, \la \A X, y \ra \, = \, \la b, y \ra = -1.
$$

$\Rightarrow:$ Suppose \eqref{p} is infeasible. We first reformulate it into a form 
described in Lemma \ref{lemma-infeas-certificate}.
Then we proceed as in the proof of the $\Rightarrow$ direction in 
Theorem \ref{theo:liftdstrong}, part \eqref{theo:liftdstrong-2} and produce a feasible solution 
of \eqref{alt-ramana}. 
\qed


\appendix 
\section{\bf Ramana's primal}
\label{section-ramana-primal} 

In this section we  study the Ramana dual of \eqref{d}, which, with some abuse of terminology we call Ramana's primal. 
The results we give here are fairly straightforward variants of results in Section \ref{section-d-ramana}, so we only sketch some of them. 
In this section we assume that \eqref{d} is feasible and the $A_i$ are linearly independent, so $\A$ is surjective. 
  
\begin{Theorem}\label{theo:main:dual}
Consider the SDP called the Ramana primal of \eqref{d} 
\begin{equation}  \label{p-ramana}  
\begin{array}{rrcl}     
	\inf        & \la C, X \ra  &     \\
	\text{s.t.} & X  & \in & \mathcal{S}_+^n + \tan( U_{n-1})   \\
	& U_0 & = & V_0  \, = \, 0                                      \\
	&  \multicolumn{3}{l}{
		\left.
		\begin{array}{rcl}
			\hspace{-5cm}  \A (U_{i} + V_{i})    & \,  = \,    & 0  \\
			\la C, U_i + V_i \ra    & \, =  \,   & 0 \\
			U_i               & \, \in \,  & \mathcal{S}_+^n \\
			V_i               & \, \in \,  & \tan( U_{i-1}) 
		\end{array}
		\right\} \quad \text{for } i =1,\dots, n-1 
	}
\end{array} \tag{P$_{\text{Ram}}$}
\end{equation}
We have that 
$$
\val \eqref{d} = \val \eqref{p-ramana},
$$
and the latter value is attained when finite. 
\qed 
\end{Theorem}  

We first need two simple propositions, both of which are  proved by elementary linear algebra.
Proposition \ref{prop-D-into-ReD} rewrites \eqref{d} in the form of \eqref{p}, i.e., in equality constrained form.
\begin{Proposition} \label{prop-D-into-ReD} 
Let $\ell = n(n+1)/2 - m, \, $ and $D_1, \dots, D_\ell$  linearly independent symmetric order $n$ matrices, such that 
$$\la A_i, D_j \ra = 0 \,\, \text{for all} \,\, i \,\, \text{and} \,\, j. $$ 
Also let  $d_j = \la D_j, C \ra$ for all $j$ and $X_0 \in \symn$ be such that 
$\A X_0 = b.$ 

Then for the SDP    
\begin{equation}\label{re-d}
\begin{array}{rl} 
	\inf  & \,\, \la X_0, Z   \ra  \\
	s.t. & \,\, \la D_j, Z \ra \, = \, \la D_j, C \ra  \, \text{for} \, j=1, \dots, \ell \\
	& \,\, Z \succeq 0, 
\end{array} \tag{\em  Re-$D$} 
\end{equation}
the following hold:
\begin{itemize}
\item[] $Z$ is feasible (optimal) in \eqref{re-d} $\Leftrightarrow$ there is $y$ feasible (optimal) in \eqref{d} such that $Z = C - \A^* y.$ 
\end{itemize}
\end{Proposition} 
\proof{} A psd matrix $Z$ is feasible in \eqref{re-d} iff $\la Z - C, D_j \ra = 0$ for all $j.$
This implies the statement about feasible solutions.

Let us next fix $Z$ and $y$ as above.
Then 
\begin{equation} \label{eqn-X0y} 
	 \la X_0, Z \ra = \la X_0, C - \A^* y \ra = \la X_0, C  \ra - \la \A X_0, y \ra = \la X_0, C \ra  - \la y, b \ra,
\end{equation}
so $ \la X_0, Z \ra  + \la y, b \ra$ is constant.  
This implies the statement about optimal solutions.

\begin{Proposition} \label{prop-AYCY-lambda} 
For $Y \in \symn$ we have 
$$
\A Y = 0, \, \la C,  Y \ra = 0 \; \Leftrightarrow \; \text{there is} \; \lambda \in \rad{\ell} \; \text{s.t.} \;  \, Y = \sum_{j=1}^\ell \lambda_j D_j \; \text{and} \; 0 = \sum_{j=1}^\ell \lambda_j \la D_j, C \ra. 
$$
\qed
\end{Proposition}
The following definition is a dual counterpart of Definition \ref{definition-RR-P}: 

\begin{Definition} \label{definition-RR-D} 
We say that \eqref{d} is in rank revealing or RR form, if the following two conditions hold:
\begin{enumerate}
\item 	 \label{definition-RR-P-Yi}  there is a regular facial reduction sequence $Y_1, \dots, Y_k$ such that 
\begin{eqnarray} \label{eqn-A*YiBi} 
	\A Y_i  & = & 0 \;\;  \text{for} \; i=1, \dots, k \\ \label{eqn-A*YB-2}  
	\la C, Y_i  \ra & = & 0 \;\;  \text{for} \; i=1, \dots, k.  
\end{eqnarray}
\item  \label{definition-RR-P-Z}  there is a slack in \eqref{d} of the form 
\begin{equation} \label{eqn-Z-Lambda} 
	\bpx 0 & 0 \\
	0 & \Lambda \epx, 
\end{equation}
where $\Lambda$ is positive definite,  and of order $n - r_{1:k}. $ Here  $r_i$ is the size of the positive definite block in $Y_i$ for $i=1, \dots, k.$ 
\end{enumerate}
We also say that the $Y_i$ in part \eqref{definition-RR-P-Yi} certify the maximum rank slack in \eqref{p}. 

\end{Definition} 
To justify the terminology of Definition \ref{definition-RR-D} we can argue just like  we did 
after Definition \ref{definition-RR-P}: if $S$ is any slack in \eqref{d}, then 
$$
\la S, Y_i \ra = 0
$$
for $i=1, \dots, k.$ Thus $Y_1$ ensures its first $r_1$ rows and columns of $S$ 
are zero;
then $Y_2$ ensures its next $r_2$ rows and columns are zero; and so on. Thus the 
$S$ slack displayed in \eqref{eqn-Z-Lambda} indeed has maximum rank. 

Lemma \ref{lemma-RR-D}   is a counterpart of Lemma \ref{lemma-RR-P}. 
Note, however, that to bring \eqref{d} into RR form we do not need elementary row operations, we only need a rotation.
\begin{Lemma} \label{lemma-RR-D}  
We can always rotate the $A_i$ and $C$ to bring \eqref{d} into RR form,  in which $k \leq n-1.$ 
\end{Lemma}
\proof{} By Lemma \ref{lemma-RR-P} we reformulate \eqref{re-d} into RR form.
Let $Q$ be the product of all rotation matrices used in the process.
First, we rotate all $D_j$ by $Q$ and also rotate all $A_i$ and $C$ by $Q, \,$ so 
$\la A_i, D_j \ra = 0$ and \eqref{eqn-X0y} still hold, and so do the conclusions of Proposition 
\ref{prop-D-into-ReD}.

Thus, after performing elementary row operations on \eqref{re-d}, we have that 
\begin{equation} \label{re-d-ref}
\begin{array}{rll} 
\inf  & \,\, \la X_0, Z   \ra  \\
s.t. & \,\, \la Y_i, Z \ra \, = \, 0  \; \text{for} \,  i = 1, \dots, k  \\
& \la D_j', Z \ra \, = \, d_j' \; \text{for} \, j \in \J \\
& \,\, Z \succeq 0, 
\end{array} 
\end{equation}
is in RR form. Here $Y_1, \dots, Y_k$ is a regular facial reduction sequence,
and the index set $\J \subseteq \{1, \dots, m \}$ indexes the other equations.  
Since the first $k$ equations in \eqref{re-d-ref} were obtained by elementary row operations, for    
$i=1, \dots, k$ we have 
\begin{equation}
Y_i = \sum_{j=1}^\ell \lambda_{ij} D_j, \; 0 = \sum_{j=1}^\ell \lambda_{ij}  \la D_j, C \ra
\end{equation} 
for some $\lambda_{ij}$ reals. 
By Proposition 	
\ref{prop-AYCY-lambda} we have $\A Y_i = 0, \, \la C, Y_i \ra = 0$ for all such $i.$
Hence \eqref{d} is in RR form, and $Y_1, \dots, Y_k$ certify the maximum rank slack in it.
\qed 

The following proposition is a counterpart of Proposition \ref{prop-reform-invariance}. 
The proof is straightforward from Proposition \ref{prop-XYQ} and the fact that elementary row operations do not change the feasible set of \eqref{p-ramana}.
\begin{Proposition} \label{prop-reform-invariance-ram-P} 
	Suppose we reformulate \eqref{p} and $Q$ is the product of all rotation matrices used in the reformulation process.
		Then 
		\begin{itemize} 
			\item 	$X$ with $\{U_i, V_i \}$ is feasible in \eqref{p-ramana} 
		before the reformulation iff $Q^\top X Q$ with 
$\{ Q^\top U_i Q, Q^\top V_i Q \}$   is feasible  after the reformulation.
	\end{itemize} 
	\qed 
\end{Proposition}

Lemma \ref{lemma-strong-primal}  is a counterpart of Lemma \ref{lemma-strong-dual}. Its 
proof is straightforward, so we omit it:

\begin{Lemma}  \label{lemma-strong-primal} 
Suppose a maximum rank slack in \eqref{d} is of the form
\begin{equation} \label{eqn-max-rank-form-dual}  
S =  Q \bpx 0 & 0 \\
0 & \Lambda \epx Q^\top, 
\end{equation}
where $Q$ is orthonormal, and $\Lambda$ is order $r$ and positive definite.
Consider the optimization problem 
\begin{equation} \label{p-strong}
\begin{split}
	\inf  & \,\, \la C, X \ra    \\
	s.t. & \,\, \la A_i, X \ra = b_i \; (i=1, \dots, m) \\ 
	& \,\,   X = Q V Q^\top \\
	& V \in \symn, \, V_{22} \in \psd{r}, 
\end{split} \tag{\mbox{${\rm P_{\rm strong}}$}}
\end{equation}
called the {\em strong primal of \eqref{d},} where $V_{22}$ stands for the lower right order $r$ block of $V.$ 
We then have 
\begin{equation} \label{eqn-val-D-val-re-P} 
\val \eqref{d} = \val \eqref{p-strong},
\end{equation} 
and $\val \eqref{p-strong}$ is attained when finite.
\qed
\end{Lemma}

Theorem \ref{theo:liftpstrong} is a counterpart of Theorem \ref{theo:liftdstrong}: it shows the feasible 
set of \eqref{p-ramana} is a lift of the feasible set of \eqref{p-strong}.
\begin{Theorem} \label{theo:liftpstrong} 
There is a $Q$ orthonormal matrix with the following properties:
\begin{enumerate}
\item \label{theo:liftpstrong-1} A maximum rank slack in \eqref{d} is of the form given in \eqref{eqn-max-rank-form-dual}.
\item \label{theo:liftpstrong-2} For any $X  \in \symn$ it holds that 
\begin{equation*}
	X \text{ is feasible in } \eqref{p-strong}
	\Leftrightarrow 
	X \text{ is feasible in } \eqref{p-ramana} \text{ with some } \{ U_i, V_i\}.
\end{equation*}
\end{enumerate}
\end{Theorem}
{\em Proof sketch}
To prove \eqref{theo:liftpstrong-1}, 
let $Q$ be the rotation constructed in Lemma \ref{lemma-RR-D}. Then after rotating \eqref{p} by this $Q$ we see 
that in the maximum rank slack in \eqref{d} the lower right block is positive definite, and the other 
elements are zero. Thus \eqref{theo:liftpstrong-1} holds. 

Next we prove \eqref{theo:liftpstrong-2}. By  
Proposition \ref{prop-reform-invariance-ram-P} we can assume 
$Q=I.$ We start with the 
$\Rightarrow$ direction. Suppose $X$ is feasible in \eqref{p-strong} and suppose 
$Y_1, \dots, Y_k$ certifies a maximum rank slack in \eqref{d}. Suppose $\Lambda_i$ is the positive definite block in $Y_i$  and let $r_i$ denote  its order for all $i. \, $ 
We decompose the $Y_i$ into $U_i + V_i$ just like we decomposed the $\A^* y^i$ in 
\eqref{eqn-define-Ui-gabor}, namely

\begin{equation} \label{eqn-define-Ui-gabor-2} 
\underbrace{	\begin{pmatrix}
	\bovermat{$r_{1:i-1}$}{\,\,\,\,\,\ti \,\,\,\,\,} & \bovermat{$\,\,\, r_{i}$}{\,\,\,\, \ti \,\,\,\,}	& \bovermat{$n - r_{1:i}$}{\,\,\,\,\,\,\, \ti \,\,\,\,\,\,\,}	\\ 
	\ti  & \, \Lambda_{i}   \,  & 0  \\ 
	\ti    &0  & 0 
	\end{pmatrix}}_{Y_i} =  
	\underbrace{	\begin{pmatrix}
	\bovermat{$r_{1:i-1}$}{\,\,\,\,\,I \,\,\,\,\,} & \bovermat{$\,\,\, r_{i}$}{\,\,\,\, 0\,\,\,\,}	& \bovermat{$n - r_{1:i}$}{\,\,\,\,\,\,\, 0 \,\,\,\,\,\,\,}	\\ 
	0  & \, \Lambda_{i}  \,  & 0  \\ 
	0    &0  & 0 
	\end{pmatrix}}_{U_i}  + \underbrace{\begin{pmatrix}
	\bovermat{$r_{1:i-1}$}{\,\,\,\,\,\ti \,\,\,\,\,} & \bovermat{$\,\,\, r_{i}$}{\,\,\,\, \ti \,\,\,\,}	& \bovermat{$n - r_{1:i}$}{\,\,\,\,\,\,\, \ti \,\,\,\,\,\,\,}	\\ 
	\ti  & \, 0  \,  & 0  \\ 
	\ti    &0  & 0 
	\end{pmatrix}}_{V_i:= Y_i - U_i}.
\end{equation}
Thus  
\begin{equation}  \label{eqn-UiVi-decomp} 
\left. \begin{array}{rcl}
\hspace{.5cm} \A (U_i + V_i)   & \,  = \,    & 0  \\
\hspace{.5cm} \la C, U_i + V_i \ra    & \,  = \,    & 0  \\
U_i               & \, \in \,  & \mathcal{S}_+^n \\
V_i               & \, \in \,  & \tan\left( U_{i-1} \right)
\end{array}
\right\} \quad \text{for } i =1,\dots, k. 
\end{equation} 
Since $k \leq n-1$ we next ``pad" the sequence $\{U_i, V_i \}$ 
with zeros. That is, we add $n-1-k$ to the index of each, and  
define 
$U_i = V_i = 0$ for $i=1, \dots, n-1-k.$ Then  \eqref{eqn-UiVi-decomp}  holds 
with $n-1$ in place of $k.$

Since $X$ is feasible in \eqref{p-strong} and $Q=I, \,$ the lower right 
order $r$ block of  $X $  is psd.  
Thus
\begin{equation} \nonumber 
X \in \psd{n} + \tan(U_{n-1}),
\end{equation}
completing the proof. 

For the $\Leftarrow$ direction, suppose $X$ with $\{U_i, V_i\}$ is feasible in 
\eqref{p-ramana}. Let $S$ be a maximum rank slack in \eqref{d}, and recall that the lower right 
order $r$ block of $S$ is positive definite, and all other elements are zero.
Then 
\begin{equation} \nonumber
\la S, U_i + V_i \ra = 0 \,\, \text{for} \,\, i=1, \dots, n-1.
\end{equation}
Thus, we argue like in \eqref{eqn-XU1XV2} and deduce that $U_i \in \psd{n, n-r}$ for all $i.$  
Hence 
$$ \tan(U_{n-1}) \subseteq \sym{n, n - r},$$ 
so the lower right order $r$ block of $X$ is psd. 
This means that  $X$ is feasible in \eqref{p-strong}, as wanted. 
\qed

\begin{Example}
	Consider again the SDP with data given in \eqref{problem-good-n=4}.
	We saw that in the dual the maximum rank slack is just the right hand side 
	\begin{equation}
		C \, = \,  \left(\begin{array}{rrrr}  
		1 & 0 & 0 & 0 \\
		0 & 1 & 0 & 0 \\
		0 & 0 & 1 & 0 \\ 
		0 & 0 & 0 & 0
	\end{array} \right).
		\end{equation}
	Thus in the strong primal we can take $Q=I$ and only the upper left $3 \times 3$ block must be psd. 
	Hence 
		\begin{equation}
		X \, = \,  \left(\begin{array}{rrrr}  
			0 & 0 & 0 & 0 \\
			0 & 0 & 0 & 1/2 \\
			0 & 0 & 0 & 0 \\ 
			0 & 1/2 & 0 & 0
		\end{array} \right).
	\end{equation}
	is feasible, and optimal in the strong primal with objective $0.$
	
	Also, in the dual
		\begin{equation}
		Y_1 \, = \,  \left(\begin{array}{rrrr}  
			0 & 0 & 0 & 0 \\
			0 & 0 & 0 & 0 \\
			0 & 0 & 0 & 0 \\ 
			0 & 0 & 0 & 1
		\end{array} \right)
	\end{equation}
	certifies the maximum rank slack. Let 
	$$ U_i = V_i = 0 \; \text{for} \;i=1,2, \; \text{and} \; U_3 = Y_1, V_3 = 0. $$
According to the proof of Theorem \ref{theo:liftpstrong}, $X$ with these $\{U_i, V_i \}$ is optimal 
in Ramana's primal of our SDP. 
  
\end{Example}
Theorem \ref{theo:liftpstrong} then directly implies Theorem \ref{theo:main:dual}. 
We can similarly translate other results in Section \ref{section-d-ramana}, e.g. Theorem \ref{thm:connect:feasibleset},  into results about \eqref{p-ramana}.
These translations are fairly straightforward to carry out, so we leave the details to the reader. 


\bibliographystyle{plain}
\bibliography{mysdp} 

\end{document}